\documentclass[preprint,12pt]{elsarticle}

\usepackage{amssymb}

\usepackage{natbib}
\usepackage{stackengine}
\usepackage{multirow}
\usepackage{hyperref}  	 
\usepackage{url}   		 
\usepackage{booktabs}  	 
\usepackage{amsfonts}  	 
\usepackage{nicefrac}  	 
\usepackage{microtype} 	 
\usepackage{lipsum}
\usepackage{fancyhdr}  	 
\usepackage{adjustbox}
\usepackage{graphicx}
\usepackage{tikz,pgfplots}
\usepackage{caption}
\usepackage{subcaption,siunitx,booktabs}
\usepackage{mathtools}
\usepackage{comment}
\usepackage{amssymb}
\usepackage{amsfonts}
\usepackage{bm}
\usepackage{amsbsy}
\usepackage{amsmath}
\usepackage{float}
\usepackage{mathtools}
\usepackage{mathrsfs}
\usepackage[linesnumbered, ruled,vlined]{algorithm2e}
\usepackage[nameinlink]{cleveref}
\usepackage{rotating}
\usepackage{adjustbox}
\usepackage{longtable}
\usepackage{pdflscape}
\usepackage[switch, modulo]{lineno}

\usepackage{empheq}
\usepackage{xfrac}
\usepackage{stmaryrd}
\usepackage{amsthm}
\usepackage{bm}

\newtheorem{theorem}{Theorem}[section]
\newtheorem{corollary}{Corollary}[section]
\newtheorem{lemma}{Lemma}[section]
\newtheorem{proposition}{Proposition}[section]
\newtheorem{definition}{Definition}[section]
\newtheorem{remark}{Remark}

\pgfplotsset{compat=1.17}

\newcommand{\bB}{\mathbb{B}}
\newcommand{\bP}{\mathbb{P}}
\newcommand{\bR}{\mathbb{R}}
\newcommand{\cB}{{\mathcal B}}
\newcommand{\cL}{{\mathcal L}}
\newcommand{\cC}{{\mathcal C}}
\newcommand{\cD}{{\mathcal{D}}}
\newcommand{\cM}{{\mathcal M}}
\newcommand{\cN}{{\mathcal N}}
\newcommand{\cP}{{\mathcal P}}
\newcommand{\cO}{{\mathcal O}}
\newcommand{\cS}{{\mathcal S}}

\newcommand{\cX}{{\mathcal X}}
\newcommand{\st}{{\text{s. t. }}}

\newcommand{\bO}[1]{{\mathcal{O}(#1)}}
\newcommand{\bsocpcomp}{BSOCP-BC\xspace}
\newcommand{\dwbc}{DW-BC\xspace}
\newcommand{\dwpwl}{DW-PWL\xspace}
\newcommand{\dwhybrid}{DW-Hybrid\xspace}
\newcommand{\dwhybrids}{DW-Hybrid*\xspace}
\DeclareMathOperator{\argmin}{argmin}

\SetKwComment{Comment}{$\triangleright$\ }{}
\SetCommentSty{}
\SetKw{Continue}{continue}

\Crefname{chapter}{Chap.}{Chaps.}
\Crefname{section}{Sect.}{Sects.}
\Crefname{proposition}{Prop.}{Props.}
\Crefname{theorem}{Thm.}{Thms.}
\Crefname{definition}{Defn.}{Defns.}
\Crefname{corollary}{Cor.}{Cors.}
\Crefname{figure}{Fig.}{Figs.}
\Crefname{observation}{Obs}{Obss.}
\Crefname{remark}{Remark.}{Remark.}

\journal{ }

\begin{document}

\begin{frontmatter}


\title{Branch and price for submodular bin packing}



\author[LIX]{Liding Xu}
\ead{lidingxu.ac@gmail.com}

\author[LIX]{Claudia D'Ambrosio}
\ead{dambrosio@lix.polytechnique.fr}

\author[LIX]{Sonia Haddad Vanier}
\ead{sonia.vanier@lix.polytechnique.fr}

\author[LIPN]{Emiliano Traversi}
\ead{traversi@lipn.univ-paris13.fr}

\affiliation[LIX]{organization={LIX CNRS, École Polytechnique, Institut Polytechnique de Paris},
   		 city={Palaiseau},
   		 postcode={91128},
   		 country={France}}

\affiliation[LIPN]{organization={LIPN CNRS, Université Sorbonne Paris Nord},
   		 city={Villetaneuse},
   		 postcode={93430},
   		 country={France}}

\begin{abstract}
The Submodular Bin Packing (SMBP) problem asks for packing unsplittable items into a minimal number of bins for which the capacity utilization function is submodular. SMBP is equivalent to chance-constrained and robust bin packing problems under various conditions. SMBP is a hard binary nonlinear programming optimization problem. In this paper, we propose a branch-and-price algorithm to solve this problem. The resulting price subproblems are submodular knapsack problems, and we propose a tailored exact branch-and-cut algorithm based on a piece-wise linear relaxation to solve them. To speed up column generation, we develop a hybrid pricing strategy to replace the exact pricing algorithm with a fast pricing heuristic. We test our algorithms on instances generated as suggested in the literature. The computational results show the efficiency of our branch-and-price algorithm and the proposed pricing techniques.
\end{abstract}



\end{frontmatter}


\section{Introduction}

Bin packing (BP) is an important combinatorial optimization problem with applications in various fields, including call centers, healthcare, container shipping, and cloud computing. These applications are typically modeled as BP problems that aim to pack unsplittable items into a minimum number of bins, with a capacity constraint on each bin.  Formally, a BP problem can be written as the following Binary Linear Programming (BIP) problem:
 \begin{subequations}
 \label{bp}
 \begin{align}
 \min \quad \sum_{j \in \cM} {y}_{{j}}, \qquad&\\
 \st \qquad \sum_{i \in \cN} \mu_i v_{ij} \le c {y}_{{j}}, \qquad& \forall j \in \cM, \label{bp.cap} \\
  \sum_{j \in \cM} v_{ij} = 1, \qquad& \forall i \in \cN, \label{bp.cover} \\
  v_{ij} \in \{0,1\}, \qquad& \forall i \in \cN, j \in \cM, \label{bp.x}\\
  y_{j} \in \{0,1\}, \qquad& \forall j \in \cM \label{bp.y},
 \end{align}
 \end{subequations}
where  $\cM:=\{1,\cdots,m\}$ is the index set of potential bins ($m$ is the number of potential bins), $\cN:=\{1,\cdots, n\}$ is the index set of items ($n$ is the number of items),  $c$ is the capacity which is the same for every bin, and $\mu_i$ is the size of item $i$.  Variable $y_j$ decides whether bin $j$ is used, and variable $v_{ij}$ indicates whether item $i$ is allocated to bin $j$.  Capacity constraints \eqref{bp.cap} stipulate that the capacities of bins are not exceeded, and set partition constraints \eqref{bp.cover} require that each item is exactly allocated to one bin.

 In many practical applications of BP, nominal item sizes $\mu$ are not revealed before the allocation decision is made, so uncertainty arises.    Probabilistic modeling of  capacity constraints \eqref{bp.cap} allows item sizes $\mu$ to be random parameters, and thus the uncertainty is taken as a probability distribution on $\mu$. We consider  two  commonly used probabilistic BP models.  The first probabilistic model is the BP with chance constraints (BPCC) \citep{song2014chance}. By assuming item sizes $\mu$ following a given (multivariate) probability distribution, BPCC requires that each capacity constraint in \eqref{bp.cap} should be respected with a probability at least $\alpha$, written as the following chance constraints \citep{charnes}:
\begin{equation}
 \label{eq.chance}
      \quad \bP(\sum_{i \in \cN} \mu_i v_{ij} \le c{y}_{{j}}) \ge \alpha, \qquad \forall j \in \cM.
	 \end{equation}
 The second probabilistic model is the distributionally robust BP (DRBP) \citep{zhang2020branch,cohen2019overcommitment}. It models the worst case of chance constraints \citep{ghaoui2003}. More specifically, given a family $ \mathscr{D}$ of probability distributions  of $\mu$, DRBP requires that each chance constraint in \eqref{eq.chance} should be respected for any probability distribution within $ \mathscr{D}$. Thereby, capacity constraints of DRBP can be formulated as the following distributionally robust constraints
 \begin{equation}
 \label{eq.chance2}
    \inf_{\mu \sim \mathscr{D}} \bP(\sum_{i \in \cN} \mu_i v_{ij} \le c{y}_{{i}}) \ge \alpha,  \qquad \forall j \in \cM.
 \end{equation}
 
 Computational optimization of BPCC and DRBP models is challenging due to probabilistic constraints.  Stochastic optimization methods can tackle mathematical optimization problems with probabilistic constraints. The sample average approximation (SAA) is a common stochastic optimization method for  chance-constrained and distributionally robust optimization problems \citep{luedtke2008sample, bertsimas2018robust}. It approximates these problems as  two/multi-stage Mixed-Integer Linear Programming (MILP) problems and computes approximate solutions that converge to an optimal solution in a probabilistic sense. Previous works, such as \cite{zhang2020branch, denton2010optimal,batun2011operating}, apply tailored SAA methods to solve BPCC and DRBP.
 
 Several recent works show that, under various assumptions on probabilistic distributions, BPCC and DRBP are equivalent to or well-approximated by a  deterministic optimization problem, namely, submodular BP (SMBP). It is shown in \cite{cohen2019overcommitment} that, BPCC has an SMBP formulation,   if item sizes $\mu$ follow independent Gaussian distributions;   SMBP also provides an upper bound for BPCC with item sizes $\mu$ under general independent distributions over bounded intervals (we note that then SMBP becomes a restriction of BPCC, and thus its solution is always feasible to BPCC.). It is shown in \cite{zhang2018ambiguous} that DRBP has an SMBP formulation, if distributions in $\mathscr{D}$ have the same mean values and the same diagonal covariance matrix.
 
 Given the applicability of the previous assumptions, SMBP is an appealing alternative formulation to BPCC and DRBP, as it can be solved optimally in a finite time, while  the convergence rate of  SAA methods for BPCC and DRBP depends on the number of samples.  SMBP already finds its applications in  cloud computing \citep{cohen2019overcommitment}, surgery planning \citep{deng2019chance}, and operating room planning \citep{wang2021chance}. The environment is highly dynamic for these applications, and uncertainty plays a significant role in practical models. These applications give rise to a need for efficient algorithms to solve SMBP.

 In this paper, we study the exact algorithms for solving SMBP. SMBP has the following Binary Nonlinear Programming formulation:

 \begin{subequations}
 \label{master.submod}
 \begin{align}
 \min \quad \sum_{j \in \cM} {y}_{{j}}, \qquad&\\
 \st \qquad \sum_{i \in \cN} a_i v_{ij} + \sigma\sqrt{\sum_{i \in \cN} b_i v_{ij} } \le c {y}_{{j}}, \qquad& \forall j \in \cM, \label{master.submod.cap} \\
  \sum_{j \in \cM} v_{ij} = 1, \qquad& \forall i \in \cN, \label{master.submod.cover} \\
  v_{ij} \in \{0,1\}, \qquad& \forall i \in \cN, j \in \cM, \label{master.submod.x}\\
  y_{j} \in \{0,1\}, \qquad& \forall j \in \cM \label{master.submod.y}
 \end{align}
 \end{subequations}
where $a_i,b_i$ are parameters inferred from the distribution of $\mu_i$. This formulation is a compact nonlinear version of \eqref{bp}. We remark that the left-hand side of the constraint \eqref{master.submod.cap} is a submodular function over $x$ \citep{atamturk2008polymatroids}, so SMBP is named after this function. A constraint in the form of \eqref{master.submod.cap} with $y_j$ fixed to 1 is called a submodular knapsack constraint.

To solve SMBP, the state-of-art exact algorithm uses general-purpose integer-programming solvers to solve its Binary Second-Order Conic Programming (BSOCP) reformulation \citep{zhang2018ambiguous}, which valid inequalities can further strengthen \cite{atamturk2008polymatroids}. The experiment  in \cite{zhang2018ambiguous} shows that small instances with item number $n$ up to $40$ and bin numbers $m$ up to 10  can be solved to optimality by this exact algorithm.

The intuition underlying this paper is that the Dantzig-Wolfe (DW) decomposition is a promising approach to tackling large-scale classical BPs: a BP is reformulated into a set cover formulation based on enumerating all feasible packing patterns; then its continuous relaxation is solved using a column generation approach \citep{gilmore1961linear}. The branch-and-price algorithm integrates column generation with the branch-and-bound algorithm. It is the state-of-the-art exact algorithm for solving DW decomposition of classical BPs \citep{wei2020new, delorme2016bin}.

We propose the first DW decomposition and set cover formulation for SMBP, and design a branch-and-price algorithm with tailored methods for solving pricing problems. After our DW decomposition of SMBP,   the nonlinearity moves to the pricing submodular knapsack subproblem, which has a linear objective function and a submodular capacity constraint. One can avoid the growing number of nonlinear constraints \eqref{master.submod.cap} in the compact formulation \eqref{master.submod}, when solving larger instances.

The DW decomposition provides a skeleton of our main algorithm. The techniques for solving general DW decomposition problems  are vast, to name a few, we refer to \cite{frangioni2013stabilized} for stabilization techniques,  \cite{coniglio2019lexicographic} for lexicographic pricing, \cite{wei2020branch} for goal cuts and early termination, and \cite{jepsen2008subset} for non-robust cuts. 
In \cite{pessoa2021solving}, a simple parameterization enables the use of several  advanced techniques in the  branch-cut-and-price VRPsolver \cite{pessoa2020generic}: automatic stabilization
by smoothing \cite{pessoa2018automation}, limited-memory rank-1 cuts \cite{pecin2017improved}, enumeration, hierarchical strong branching over accumulated resources \cite{gelinas1995new}, and
limited discrepancy search diving heuristics.
In this paper, we focus on algorithmic innovation that exploits the specific nonlinear structure of pricing problems and new techniques to speed  up the convergence of column generation.

 As the study in \cite{zhang2018ambiguous} for the compact formulation, the nonlinearity is a crucial feature for model representability, so it is unavoidable and needs a special algorithmic treatment. In our case, pricing problems have submodular knapsack constraints involving nonlinear functions. We give two different views of nonlinearity. First, we can represent the submodular knapsack constraint via second-order constraints, which many general solvers then accept. Alternatively, we propose a non-convex Mixed-Binary Quadratically Constrained Programming (MBQCP) formulation for the submodular knapsack. A Piece-Wise Linear (PWL)  function is linear in each partition of its domain and can be modeled by a MILP formulation \citep{vielma2010}. PWL functions have been used to approximate or relax non-convex Mixed-Integer Nonlinear Programming (MINLP) problems \citep{geissler2012using}.  Despite its non-convexity, the critical feature of the MBQCP formulation is that its only nonlinear function is a univariate quadratic function, which is easy to approximate using a PWL function. We construct the PWL relaxation for the submodular knapsack and combine it with cutting planes to form an exact PWL relaxation-based branch-and-cut (PWL-B$\&$C) algorithm.

 The submodular knapsack  is essential as it models the chance-constrained knapsack problem \citep{goyal2010ptas}. We thus provide an approach for solving submodular knapsack problems different  from the pure valid inequality approach in \cite{atamturk2008polymatroids,atamturk2009submodular}.
 
 We propose  several strategies to accelerate the convergence of the branch-and-price algorithm, i.e., improve primal and dual bounds. The Farley bound \cite{farley1990note,vance1994solving}  is an early valid dual bound before the termination of the column generation procedure \cite{wei2020new, gleixner2020price}.
The formula for the Farley bound imposes a condition on whether an exact pricing algorithm can improve the current dual bound. If the condition is not satisfied, the exact pricing algorithm is unnecessary, so we can use fast pricing heuristic. Our branch-and-price algorithms use a hybrid pricing strategy to speed up the column generation procedure. The hybrid pricing strategy is thus an intermediate between exact pricing and heuristic pricing strategies \cite{blanco2023branch}.
 
 There are few publicly available instances of SMBP problems. In \cite{cohen2019overcommitment}, there is a method to generate instances from BPCC and DRBP under various distributions. We generate instances of three different scales  by this method and conduct computational experiments on them. We implement our branch-and-price algorithm for the set cover formulation of SMBP and find that it outperforms existing methods, which solve the compact BSOCP formulation \cite{zhang2018ambiguous}. Our core innovation,  PWL-B$\&$C pricing algorithm, and hybrid pricing strategy, significantly improve the branch-and-price algorithm.

\subsection{Literature review}
As mentioned, there are several steps of transformation from BP with uncertainty to SMBP. We review these transformations and algorithms for solving associated transformed models.

The surgery planning problem is a typical application of BP with uncertainty in healthcare, where the surgery duration (item size) is assumed to be stochastic. Some pioneering works
 \citep{ denton2010optimal, batun2011operating} allow violations to capacity constraints (the left-hand side of \eqref{bp.cap} thus can be greater than the capacity) rather than consider chance constraints.  To minimize these violations, they  use a penalty approach by adding the expectation of the sum of these violations into the objective function. Therefore, the transformed BP model is a standard stochastic optimization problem called stochastic BP (SBP). SBP can be further modeled and solved as a stochastic two-stage mixed-integer programming problem: the first stage variables are the bin variables $y$, the second stage variables are the item variables $v$, and the expected violation is the second stage objective.  In some works \citep{cardoen2010operating, deng2020}, only the expected penalty is considered in the models.

 Compared to SBP, BPCC, and DRBP can control the violation of each capacity constraint with  a guaranteed probability bound, and thus they are more accurate models for BP with uncertainty. To solve BPCC and DRBP, there are approximation algorithms  and exact algorithms. (Sampling-based) approximation algorithms usually converge  asymptotically to an optimal solution when the sampling number increases (as SAA methods), and exact algorithms usually converge in  finite time. However, exact algorithms are mostly available for deterministic optimization problems.
 
 In \cite{shylo2013}, a variant of  BPCC in surgery planning is studied:   items (surgeries) are allocated to  a  given set of bins (time blocks), the goal is to minimize the sum of expected capacity residuals (undertime), subject to chance constraints for overuse of  bins' capacities (overtime). Assuming that the operation duration follows a multivariate normal distribution, the authors reformulate the problem as a deterministic optimization problem containing a convex objective and submodular capacity constraints. In \cite{song2014chance}, a special BPCC with the probability distribution over finite support is studied. It has  a BIP formulation, which an exact algorithm can solve.  As mentioned before, for various probabilistic distributions \cite{cohen2019overcommitment,zhang2018ambiguous},  BPCC and DRBP admit a deterministic SMBP reformulation, which can be solved by exact algorithms.

Regarding solution algorithms, an SAA-based   algorithm \cite{zhang2020branch} can solve BPCC approximately, whose scenario subproblem is solved exactly by a DW decomposition approach. In addition, SMBP reformulation of DRBP  \cite{zhang2020branch} can be approximated by a MILP, which is solved by a DW decomposition approach. The only tailored exact algorithm for BPCC and DRBP  \cite{zhang2018ambiguous} solves their compact SMBP reformulations. As for the deterministic variant of BPCC, the authors of   \cite{shylo2013} propose an exact outer approximation algorithm enhanced with PWL relaxation of submodular knapsack constraints, and the algorithm is a multi-search tree method, i.e., an underlying MILP solver will be called multiple times. In conclusion, no exact algorithm based on DW decomposition exists to solve SMBP. Meanwhile, DW decomposition is  already used for various approximated problems.

This exception  may be due to the lack of efficient exact algorithms to solve submodular knapsack problems.
We note that in \cite{zhang2020branch}, DW decomposition is \textit{de facto} applied to a MILP problem, and thus the  pricing problems are also MILPs such as classical knapsack problems \citep{cacchiani2022knapsack}, which can be solved efficiently by general-purpose integer programming solvers.  The efficiency of general solvers is mostly due to the lifted cover inequalities, which are strong valid inequalities for knapsack polytope and can be constructed via sequence-independent lifting \citep{gu1999lifted}. However, for submodular knapsack, the computation of lifted cover inequalities is not tractable  \citep{atamturk2009submodular}. As we know from the literature,  the tailored algorithm can be much better than general solvers for many variants of classical knapsack, because these algorithms can exploit more problem structures than general solvers. To name a few, we refer  to the quadratic knapsack \citep{caprara1999, furini2019theoretical}, the multidimensional knapsack \citep{Puchinger2010}, and the quadratic multi-knapsack \citep{Bergman19, olivir2021}.

Although general solvers are almost as complex as a black box for users, we can at least understand how they solve the submodular knapsack. The submodular knapsack can be reformulated as a BSOCP problem, which is an acceptable formulation to  \texttt{CPLEX} \citep{bliek1u2014solving} and \texttt{SCIP} \citep{berthold2012extending}. These solvers implement  LP outer approximation-based branch-and-cut (LP-B$\&$C) algorithm \citep{coey2018outer} to solve the BSOCP or general Mixed-Integer Second-Order Conic Programming (MISOCP) problems. The LP outer approximation is sometimes called the polyhedral outer approximation (or polyhedral relaxation).  In fact, any second-order conic program (SOCP) is polynomially reducible to a linear program \cite{ben2001polyhedral}. As for a submodular knapsack constraint, general solvers will linearize it into an intersection of a set of classical knapsack constraints, thus, inefficiency arises if too many linearizations are applied.

On the other hand, there are alternative exact approaches to solve some classes of nonconvex MINLPs using PWL relaxations  \citep{d2012algorithmic, geissler2012using}. For example, \cite{d2012algorithmic} obtains a convex MINLP relaxation for nonconvex MINLPs with separable nonconvex functions. The authors distinguish between convex and concave parts and then convexify the concave parts by PWL functions.

Regarding the submodular knapsack, its BSOCP (a convex MINLP) formulation has a nonlinear function over all the problem variables, which is difficult to approximate when  the dimension is high. On the other hand,  the nonconvex MBQCP formulation, where the only nonlinear function is a univariate quadratic function on a slack variable. The resulting PWL relaxation in the experiment is stronger than pure polyhedral relaxation. In our case, we will show that the quadratic function can be approximated in a ``dimension-free" way, since the nonlinearity is concentrated on a single variable. \cite{shylo2013} only uses  PWL relaxations to approximate the submodular knapsack. It requires refining PWL relaxations to achieve convergence, so their multi-search tree algorithm needs to restart the MILP solver from scratch in each iteration. In contrast, we prefix PWL relaxations and use cutting planes to achieve convergence. Therefore, our  single-search tree does not need to restart the MILP solver.

 We look at the recent development of DW decomposition and  branch-and-price algorithm for solving MINLPs
 (see \cite{allman2021branch}), such as recursive circle packing (RCP) problems  \citep{gleixner2020price}, binary quadratic problems \citep{ceselli2022dantzig}, and facility location with general nonlinear facility cost functions \citep{ni2021branch}.  There may be several ways to divide a MINLP into master and subproblems, so a MINLP may admit different DW decompositions. In \cite{ceselli2022dantzig}, the authors study the strengths of different DW decompositions for binary quadratic problems. In most cases, after applying the DW decomposition to the compact MINLP formulation, the master problem is a MILP, and the pricing problems are MINLPs.  Since pricing problems are solved in thousands of iterations, \cite{gleixner2020price} shows that any improvement in the pricing algorithm can speed up the convergence of column generation.

\subsection{Contribution}

In summary, our contribution in this paper is threefold. As far as we know, the previous work applies DW decomposition for an approximated MILP for SMBP, and thus nonlinearity is not considered in the solving process. So, we are the first to apply  DW decomposition for SMBP with the nonlinearity considered. Built on the basic DW decomposition and branch-and-price algorithm, we develop a new hybrid pricing strategy technique to speed up the column generation, which can avoid computationally expensive exact pricing while not worsening the dual bound.  Second, for pricing submodular knapsack problems, we  propose a new MBQCP formulation and its PWL relaxation, and design a new PWL-B$\&$C algorithm  as an alternative  exact algorithm to the conventional LP-B$\&$C algorithm, which is  based on valid inequalities. Finally, we perform computational experiments on many instances to evaluate the proposed algorithms. The computational results show that our tailored branch-and-price algorithms for DW reformulation outperform the conventional branch-and-cut algorithm for BSOCP formulation implemented in a state-of-art commercial solver; and the PWL-B$\&$C algorithm can be a standalone algorithm for submodular knapsack. The source code and benchmark are released on our project website \href{https://github.com/lidingxu/cbp}{https://github.com/lidingxu/cbp}.

\subsection{Outline}
 This paper is organized as follows. In \Cref{sec.prob}, we describe the set cover formulation of SMBP. In \Cref{sec.bp}, we introduce the critical components of our branch-and-price algorithm: the branching rule, column generation, dual bound computation, initial columns, and primal heuristics. In \Cref{sec.price}, focusing on solving the pricing problem, we present the pricing heuristic, reformulations of the pricing problem, PWL relaxation, the exact pricing algorithm, and the hybrid pricing strategy. In \Cref{sec.experi}, we show the computational results of the proposed algorithms for instances generated from the literature and analyze their performance. In \Cref{sec.conc}, we end this paper with a conclusion and future research directions.

\section{Set cover formulation}
\label{sec.prob}
In this section, we propose a new set cover formulation for SMBP. The formulation is derived similarly to the DW decomposition of the classical linear BP \citep{delorme2016bin}. This formulation can be solved efficiently by a branch-and-price algorithm.

A column \(p\) is defined by a binary vector as \((d_{1p}, d_{2p}, \hdots,d_{np})\), where \(d_{ip} = 1\) if item \(i\) is contained in the column \(p\). A column is called \textit{feasible} if the combination of its items can fit into a bin, i.e., satisfies the submodular capacity constraint \eqref{master.submod.cap}. The set cover formulation is based on enumerating all feasible columns, the number of which can be exponential to the number of items.

\textbf{Set notation:}

\begin{itemize}
 \item \(\cP\): the set of all feasible columns.
\end{itemize}

\textbf{Decision variables:}
\begin{itemize}
 \item $\lambda_p = \begin{dcases} 1, & \textup{if column \(p\) is used by the solution}\\
  0, & \textup{otherwise}
 \end{dcases}$ for $p \in \cP$.
\end{itemize}

We obtain the following set cover formulation for SMBP:

 \begin{subequations}
 \label{master.dw}
 \begin{align}
 \min \quad \sum_{p \in \cP} \lambda_p, \qquad&\\
 \st \qquad \sum_{p \in \cP} d_{ip} \lambda_p \ge 1, \qquad& \forall i \in \cN, \label{master.dw.sc}\\
 \qquad \lambda_p \in \{0,1\}, \qquad & \forall p \in \cP.
 \end{align}
 \end{subequations}
The set cover constraint \eqref{master.dw.sc} specifies that each item \(i\) ($i \in \cN$) is contained in at least one bin. The set cover reformulation already finds applications in vehicle routing \cite{pecin2017improved} and unsplittable multi-commodity flows \cite{xu2022branch} problems.

The compact formulation \eqref{master.submod} is a MINLP, but the set cover formulation \eqref{master.dw} is a MILP. Moreover, the number of nonlinear constraints in the compact formulations equals the number of potential bins. The nonlinearity of the set cover formulation is \textit{de facto} `hidden' in the pricing subproblems, and each pricing subproblem has only one nonlinear constraint.

\begin{remark}
\label{rm.bsocp}
(Modelling BSOCP constraints)
We give a way to obtain a BSOCP formulation of the constraint $\sum_{i \in \cN} a_i v_{ij} + \sqrt{\sum_{i \in \cN} b_i v_{ij} } \le d$, where $d = c y_i$ or $d = c$. Since $v_{ij} \in \{0,1\}$, the square root $\sqrt{\sum_{i \in \cN} b_i v_{ij} }$ equals  $\sqrt{\sum_{i \in \cN} b_i v^2_{ij} }$. Let $d' := d - \sum_{i \in \cN} a_i v_{ij}$, then the constraint is equivalent to a SOCP constraint $ \sqrt{\sum_{i \in \cN} b_i v^2_{ij} } \le d'$. One can transform the SOCP constraint into several 3d SOCP constraints like $x_1 \ge \sqrt{x_2^2 +x_3^2}$, which is acceptable by \texttt{CPLEX} \cite{bonami2015recent}.
\end{remark}

Through the above discussion, we can obtain a BSOCP formulation of SMBP \eqref{master.submod}.
When comparing two formulations, a formulation is said to be ''stronger'' if it gives a better dual bound.

\begin{proposition}
\label{prop.strong}
The linear relaxation of the set cover formulation \eqref{master.dw} is stronger than the continuous SOCP relaxation of the BSOCP formulation of \eqref{master.submod}.

\end{proposition}
The proof can be found in  \ref{proof.prop.strong}.

\section{Branch and price}

\label{sec.bp}
Solving the set cover formulation with an exponential number of binary variables is challenging.
In this section, we present an exact branch-and-price algorithm for solving the set-cover formulation of SMBP. The branch-and-price algorithm integrates column generation with the branch-and-bound algorithm to efficiently solve the LP relaxation. In the following subsections, we describe the important steps of our branch-and-price algorithm: the branch rule, the column generation, the primary heuristic, and the dual bound computation.

\subsection{Branching rule}

 Our branch-and-price algorithm uses the Ryan/Foster branching rule \citep{ryan1981integer}. The branching rule selects a pair of items $i_1 \in \cN$ and $i_2 \in \cN$ that must either be packed together or not packed together. We denote by

\begin{itemize}

 \item $\cS$: the set of item pairs that are forced to be packed together such that, if a column $p$ respects $\cS$, then for $(i_1, i_2) \in \cS$, $d_{i_1p} = d_{i_2p}$;

 \item $\cD$: the set of item pairs that are not allowed to be packed together such that, if a column $p$ respects $\cD$, then for $(i_1, i_2) \in \cD$, $d_{i_1p} + d_{i_2p} \le 1$.

\end{itemize}

Indeed, $(\cS, \cD)$ exactly describes the branching decisions made for each node of the search tree, whose  nodes are constructed and selected by \texttt{SCIP}'s internal rules \cite{achterberg2008constraint} in our implementation.
We denote by
\begin{equation*}
 \cP_{\cS, \cD} := \{p \in \cP \,|\, \forall (i_1, i_2) \in \cS \; d_{i_1p} = d_{i_2p} \land \forall (i_1, i_2) \in \cD\; d_{i_1p} + d_{i_2p} \le 1\}
\end{equation*}
the set of feasible columns respecting branching constraints induced by $(\cS, \cD)$. We refer to $\cP_{\cS, \cD}$ as the $(\cS, \cD)$-feasible columns.

At each node of the search tree, the set cover problem \eqref{master.dw} is restricted to the branching decision set \((\cS, \cD)\), i.e., it follows as

\begin{subequations}
 \label{master.dwsd}
 \begin{align}
 \min \quad \sum_{p \in \cP_{\cS, \cD} } \lambda_p, \qquad&\\
 \st \qquad \sum_{p \in \cP_{\cS, \cD} } d_{ip} \lambda_p \ge 1, \qquad& \forall i \in \cN, \label{master.dwsd.sc}\\
 \qquad \lambda_p \in \{0,1\}, \qquad& \forall p \in \cP_{\cS, \cD} .
 \end{align}
 \end{subequations}

The above problem \eqref{master.dwsd} is called the \textit{master problem}, and its LP relaxation is called the \textit{master LP problem}.

Given a solution $\lambda$ of the LP relaxation, if $\lambda$ is not integral, the branching rule chooses an item pair to branch. It first creates an $n$-by-$n$ matrix, and computes its entries as $M_{i_1i_2}= \sum_{p \in \cP_{\cS,\cD}:  d_{i_1p} = d_{i_2p} = 1}\lambda_p$ for all $i_1,i_2 \in \cN$.  Since, for an integral solution, $M_{i_1i_2}$ must be either 0 or 1, the branching rule chooses the most fractional entry $(i'_1,i'_2)$ such that $i'_1,i'_2 = \argmin_{i'_1, i'_2 \in \cN} |0.5-M_{i'_1i'_2}|$. Then, the rule adds $(i'_1,i'_2)$ to $\cS,\cD$, respectively.

\subsection{Column generation}

\label{sec.col}

We present a column generation method to solve the master LP problem.

The column generation procedure starts with a subset of \((\cS, \cD)\)-feasible columns of the master LP problem, adds columns, and solves the restricted LP iteratively. Given a subset \(\cP'_{\cS,\cD}\) of \(\cP_{\cS,\cD}\), the corresponding restricted LP problem, namely the \textit{Restricted Master LP } (RMLP) problem, is

 \begin{subequations}
 \label{master.dw.rel}
 \begin{align}
 \min \quad \sum_{p \in \cP'_{\cS,\cD}} \lambda_p, \qquad&\\
 \st \qquad \sum_{p \in \cP'_{\cS,\cD}} d_{ip} \lambda_p \ge 1, \qquad& \forall i \in \cN, \label{master.dwrel.sc}\\
 \qquad \lambda_p \ge 0, \qquad& \forall p \in \cP'_{\cS,\cD}.
 \end{align}
 \end{subequations}

After solving the RMLP, let \(\pi_i\) be the dual variable associated with the \(i\)-th constraint \eqref{master.dwrel.sc}. The reduced cost for a column $p \in \cP_{\cS,\cD}$ is \(r_p := 1 - \sum_{i \in \cN} \pi_i d_{ip}\). If there is a column \(p \in \cP_{\cS,\cD} \setminus \cP'_{\cS,\cD}\) whose reduced cost \(r_p\) is negative, then adding \(p\) to \(\cP'_{\cS,\cD}\) could reduce the objective value of the RMLP. Otherwise, the solution for the RMLP is also optimal for the master LP problem. The column with the most negative reduced cost is determined by solving a pricing problem.

Before the column generation procedure is applied to the current node, the items that can only be packed together are combined into the set $\cS$ using a preprocessing process. Let the new item set be $\cN'$, $a',b'$ be the merged parameters, and the new conflict relation be $\cD'$. Preprocessing leads to a smaller pricing problem, which can be formulated to a \textit{submodular knapsack problem with conflicts}:
 \begin{subequations}
 \label{pricing.ref}
 \begin{align}
  \max \quad \sum_{i\in \cN'} \pi'_i x_i, \qquad& \\
  \st \qquad \sum_{i \in \cN'} a'_i x_i + \sigma\sqrt{\sum_{i\in \cN'} b'_i x_i } \le c,\qquad& \label{pricing.ref.cap}\\
  x_{i_1} + x_{i_2} \le 1, \qquad & \forall (i_1, i_2) \in \cD',\label{pricing.ref.conflict}\\
  x_i \in \{0,1\},\qquad & \forall i \in \cN'.
 \end{align}
 \end{subequations}

If the optimal value $\sum_{i\in \cN'} \pi'_i x_i > 1$, then the corresponding column has a negative reduced cost $1 - \sum_{i\in \cN'} \pi'_i x_i $ and is added to the RMLP. Otherwise, the solution of the RMLP is optimal for the master LP, and the current node is solved. The details of the pricing algorithms can be found in \Cref{sec.price}.

Since pricing problems may not be solved optimally within a time limit, a pricing algorithm may find an existing column in $\cP'_{\cS,\cD}$ or a column with positive reduced cost. Therefore, adding the column does not improve the RMLP and the column generation procedure is aborted. The following simple constraint can exclude existing solutions from the pricing problem, thus reducing the search space:
 \begin{equation}
     \label{cons.exclude}
     \sum_{i\in \cN'} \pi'_i x_i \ge  1 + \epsilon,
 \end{equation}
 where $\epsilon$ is a sufficiently small positive real number. Exact algorithms can easily add this constraint to exclude existing columns in $\cP'_{\cS,\cD}$. This constraint  also guarantees that if  solutions of negative reduced costs exist, then exact algorithms can find one of them.

\subsection{Primal heuristics}
\label{sec.pheur}
We discuss primal heuristics that help find primal feasible solutions for the set covering formulation. We use two heuristics: The first heuristic uses an approximation algorithm to find a primal solution that forms a set $\cP'$ of initial columns, and the second heuristic attempts to find a primal solution once a column has been generated and added to $\cP'$.

 \cite{cohen2019overcommitment} propose approximation algorithms to find a feasible solution with $8/3$-ratio to the optimal solution to the submodular bin packing. Their algorithms are greedy and easy to implement, so we employ these algorithms as the first heuristic.

During column generation, each generated column could be combined with the previous columns in $\cP'$ into a  primal feasible solution. Our second primal  heuristic is similar to the greedy column selection heuristic in \cite{lubbecke2012primal,joncour2010column}.  Once a column is generated, we force it into a potential solution. Then, we greedily select an existing column from $\cP'$ that packs the maximum number of unpacked items until all items are packed.  We note that the heuristic may find columns that do not improve the RMLP.

\subsection{Dual bound computation}

For an optimization problem, a dual bound certifies the optimality of a solution. In the branch-and-price setting, a local dual bound at each node of the search tree is a lower bound on the optimum of the master problem \eqref{master.dwsd}. The algorithm uses the local dual bound to fathom the node or select branch nodes.

The optimum of the master LP problem is a local dual bound. However, the column generation procedure usually needs to solve many pricing problems to converge to this optimum. At each iteration of the column generation procedure, another local dual bound is available. This bound is referred to in the literature as \textit{Farley bound}. The following lemma illustrates how this bound can be computed.

\begin{lemma}[ \cite{farley1990note,vance1994solving}]

\label{lem.valildlbd}

Let \(v_{\mathrm{MP}}\) be the optimum of the master LP, let \(v_{\mathrm{RMLP}}\) be the optimum of the RMLP, let \(v_{\mathrm{price}}\) be a dual bound for the pricing problem \eqref{pricing.ref}, and let \(v_{\mathrm{F}} \coloneqq \frac{v_{\mathrm{RMLP}}}{v_{\mathrm{price}}}\) be the Farley bound. Then, $ v_{\mathrm{F}} \le v_{\mathrm{MP}}$, and thus $v_{\mathrm{F}}$ is a local dual bound.

\end{lemma}

The computation of the Farley bound requires a dual bound on the pricing problem, obtained using an exact pricing algorithm. The branch-and-price algorithm holds a local lower bound $v_{\mathrm{ld}}$ at each search tree node. After solving each pricing problem, the branch-and-price algorithm updates $v_{\mathrm{ld}}$ according to the following rule:

\begin{equation*}
 \label{eq.update}
 v_{\mathrm{ld}} = \max \{ v_{\mathrm{F}}, v_{\mathrm{ld}}\}.
\end{equation*}

Early termination rules of \cite{wei2020new} can compare the local dual bound and the primal bound to improve the branch-and-price algorithm. The rules exploit integrality and can terminate column generation earlier than the classical algorithm. We implement these rules in our branch-and-price solver. 

\section{Solving the pricing problem}
\label{sec.price}
In this section, we present solution methods for the pricing problem. The proposed algorithms can be implemented as a stand-alone solver for the submodular knapsack problem.

We first present a fast pricing heuristic. We then present two formulations of the submodular knapsack problem (with conflicts): a convex BSOCP formulation and a non-convex MBQCP formulation. The convex BSOCP formulation is solved in our experiments for a comparative study. The PWL method is a way to approximate nonlinear functions (or relax under some conditions) by linear functions in its subdomain. We derive a PWL relaxation of the MBQCP formulation and develop an exact PWL-based branch-and-cut algorithm (PWL-B$\&$C) for the pricing problem.

To speed up column generation, we also present a hybrid pricing strategy  that can replace the exact pricing algorithm with a fast pricing heuristic.

\subsection{Pricing heuristic}
\label{sec.heur}

We propose a fast heuristic, the fixing-greedy heuristic. This heuristic is used by the hybrid pricing strategy to speed up the column generation procedure.

The fixing-greedy heuristic is based on the best-fit-greedy algorithm. The best-fit-greedy algorithm adds an item per iteration only if it does not conflict with the previously added items, as long as the capacity is not exceeded. The heuristic keeps
\begin{itemize}
    \item $\Delta$:   the set of items added to the bin, which is initially empty.
\end{itemize}

 At each iteration, the best-fit greedy heuristic has the following steps:

\begin{enumerate}
    \item computes the sum of \(a'_i\) and  the sum of
\(b'_i\) of added items, i.e.,  \(A := \sum_{i \in \Delta} a'_i\) and \(B := \sum_{i \in \Delta} b'_i\);
\item find the set $\overline{\Delta}:=\{i \in \cN' \setminus \Delta: A + a'_i + \sigma \sqrt{B + b'_i} \le c\}$ of items that can be added to the bin;
\item if $\overline{\Delta} = \emptyset$, exits and outputs $\Delta$;
\item for each unadded item \(i \in \overline{\Delta}\), computes  the incremental capacity usage \(\gamma_i :=(A + a'_i + \sigma \sqrt{B + b'_i}) -(A + \sigma \sqrt{B})\), and the profit-over-usage ratio \(r_i := \frac{\pi'_i}{\gamma_i}\);
\item adds the unadded item with the maximum \(r_i\) into $\Delta$.
\end{enumerate}

The fixing-greedy heuristic  enforces, for each time, an item
in \(\cN'\) to be in the solution, runs the best-fit greedy algorithm,  and outputs the best solution.

\subsection{BSOCP formulation}

The Binary Second-Order Conic Programming formulation of the pricing problem \eqref{pricing.ref} is similar to the BSOCP formulation of SMBP \eqref{master.submod}.

Applying the same technique in \Cref{rm.bsocp},
the BSOCP formulation of the pricing problem is:
    \begin{subequations}
    \label{pricing.ref.bsocp}
    \begin{align}
  		 \max \quad \sum_{i\in \cN'} \pi'_i x_i, \qquad& \\
   	 \st \qquad  \sum_{i \in \cN'} a'_i x_i  + \sigma\sqrt{\sum_{i\in \cN'} b'_i x_i^2 } \le c,\qquad& \label{pricing.ref.bsocp.cap}\\
   	 x_{i_1} + x_{i_2} \le 1, \qquad & \forall (i_1, i_2) \in \cD', \label{pricing.ref.bsocp.conflict}\\
   	 x_i \in \{0,1\},\qquad & \forall  i  \in \cN'.
    \end{align}
    \end{subequations}
Where \eqref{pricing.ref.bsocp.cap} can be represented by 3d second-order conic constraints. The BSOCP formulation \eqref{pricing.ref.bsocp} is a convex MINLP formulation.

In this section, we analyze the polyhedral outer approximation of the BSOCP formulation \eqref{pricing.ref.bsocp}  and show that a finite number of cutting planes is sufficient to define an exact MILP reformulation of the BSOCP formulation \eqref{pricing.ref.bsocp}.

To simplify the presentation, we use the following notation:
\begin{itemize}
    \item the left-hand side of \eqref{pricing.ref.bsocp.cap}: \begin{equation*}
\label{eq.subfunc}
    f(x) := \sum_{i \in \cN'} a'_i x_i  + \sigma\sqrt{\sum_{i\in \cN'} b'_i x_i^2 };
\end{equation*}
\item the binary set defined by \eqref{pricing.ref.bsocp.cap}: \begin{equation*}
\label{eq.cC}
    \cC := \{x \in \{0,1\}^{\cN'}:  f(x) \le c\};
\end{equation*}
\item the continuous relaxation of $\cC$: \begin{equation*}
\label{eq.cCrel}
    \overline{\cC} := \{x \in [0,1]^{\cN'}:  f(x) \le c\}.
\end{equation*}

\end{itemize}
 
Since $f$ is convex,  $ \overline{\cC}$ is convex. We also note that the convex hull of $\cC$ is a polytope. A set $\cO$ is a polyhedral outer approximation of $\cC$, if $\cO$ is a  polyhedron and $\cC \subset \cO$. A  polyhedral outer approximation can be constructed as follows. Define a  linearization of  $f$ at some $\hat{x}$ in the domain of $f$ by $ \cL^{f}_{\hat{x}}(x) \coloneqq f(\hat{x}) + \nabla f(\hat{x}) ^\top (x - \hat{x})$.  Since $f$ is convex,  $ \cL^{f}_{\hat{x}}$ is an under-estimator of $f$, i.e., $ \cL^{f}_{\hat{x}}(x) \le f(x)$ for any $x$. Hence,  $\cL^{f}_{\hat{x}}(x) \le c$ is a linear inequality valid for $f(x) \le c$.

A polyhedral outer  approximation $\cO$ is said  \textit{exact}, if $\cO \cap \{0,1\}^n = \cC$. So, solving the optimization problem over an exact polyhedral outer approximation  with binary and conflict constraints is equivalent to solving the submodular knapsack problem with conflicts. Next, we identify a family of valid inequalities that give an exact polyhedral outer  approximation. Each of these valid inequalities corresponds to a binary point not in $\cC$.

\begin{theorem}
\label{lem.sep}
Given a point \(\hat{x} \in \{0,1\}^{\cN'} \),  the following  inequality is valid for \(\cC\) and $\overline{\cC}$:
\begin{equation}
\label{eq.outercut}
	 \sum_{i \in \cN'} a'_i x_i  + \frac{\sigma}{  \sqrt{\sum_{i\in \cN'} b'_i \hat{x}_i}} {\sum_{i\in \cN'} b'_i \hat{x}_i x_i}\le c.
\end{equation}
Let $$\cO = \{x \in [0,1]^{\cN'} :  \sum_{i \in \cN'} a'_i x_i  + \frac{\sigma}{  \sqrt{\sum_{i\in \cN'} b'_i \hat{x}_i}} {\sum_{i\in \cN'} b'_i \hat{x}_i x_i}\le c, \; \forall \hat{x} \in \{0,1\}^{\cN'} \setminus \cC \}.$$
Moreover,
\begin{enumerate}
    \item if $\hat{x} \notin \cC$, the valid inequality is violated by $\hat{x}$;
    \item  $\cO$ is exact, and $\cC = \cO \cap \{0,1\}^{\cN'}$.
\end{enumerate}
\end{theorem}
\begin{proof}

Since function  $f$  is  convex,  it follows that
\begin{equation*}
    \cL^f_{\hat{x}}(x)   \le f(x) \le c.
\end{equation*}

Moreover,
\begin{equation*}
\begin{split}
   &   \cL^f_{\hat{x}}(x)\\
   = &f(\hat{x}) + \nabla{f(\hat{x})} ^\top (x -\hat{x})\\
    =&  \sum_{i \in \cN'} a'_i x_i   +  \sigma\sqrt{{\sum_{i\in \cN'} b'_i \hat{x}_i^2}}+  \frac {\sigma}{ \sqrt{\sum_{i \in \cN'} b'_i \hat{x}_i^2}}\sum_{i \in \cN'}  b'_i \hat{x}_i (x_i-\hat{x}_i)\\
    =& \sum_{i \in \cN'} a'_i x_i  + \sigma\sqrt{{\sum_{i\in \cN'} b'_i \hat{x}_i^2}}+ \frac{\sigma}{  \sqrt{\sum_{i\in \cN'} b'_i \hat{x}_i^2}} {\sum_{i\in \cN'} b'_i \hat{x}_i x_i} -  \frac{\sigma}{  \sqrt{\sum_{i\in \cN'} b'_i \hat{x}_i^2}} {\sum_{i\in \cN'} b'_i \hat{x}_i \hat{x}_i}\\
    =&  \sum_{i \in \cN'} a'_i x_i+ \frac{\sigma}{  \sqrt{\sum_{i\in \cN'} b'_i \hat{x}_i}} {\sum_{i\in \cN'} b'_i \hat{x}_i x_i}
\end{split}
\end{equation*}
where the last equation follows from the fact that $\hat{x}$ is binary.
 
Therefore, inequality \eqref{eq.outercut} in the statement is valid for $\cC$. The left-hand side  of  inequality \eqref{eq.outercut} evaluated at $\hat{x}$ is $ \sum_{i \in \cN'} a'_i \hat{x}_i + \sigma \sqrt{ \sum_{i\in \cN'} b'_i \hat{x}_i} $ which is by hypothesis is at least $c$, so  $\hat{x}$ violates the inequality.

Let us consider $x^\ast \in \{0,1\}^{\cN'}$. If $x^\ast\notin \cC$, then  $x^\ast$ violates the $\cL^f_{x^\ast}(x)   \le c$ which is a facet defining inequality of $\cO$, then  $x^\ast \notin \cO$. Hence, $x^\ast \in \cO$ implies that $x^\ast \in \cC$. If $x^\ast \in  \cC$, since $\cO$ is a polyhedral outer approximation of $\cC$, $x^\ast$ must be in $\cO$. Therefore,  $\cC = \cO \cap  \{0,1\}^{\cN'}$.  \end{proof}

Looking at the above theorem, we find that each binary point not in $\cC$ gives rise to a valid inequality separating it from $\cC$. Moreover, binary points in $\cC$ satisfy these valid inequalities, i.e., they are in the polyhedral outer approximation $\cO$. We define two sets related to the polyhedral outer approximation $\cO$.  The \textit{generating set} is defined as
\begin{equation}
\label{eq.gen}
    \cX:=\{\hat{x} \in \{0,1\}^{\cN'}: \hat{x} \notin \cC\},
\end{equation} because it generates the following \textit{cut coefficient set}:
\begin{equation}
    \label{corecuts}
    \Theta := \left\{\theta \in \bR^{\cN'}: \exists \hat{x}  \in \cX \,  \forall i \in \cN' \, \theta_i=  a'_i  + \frac{\sigma}{  \sqrt{\sum_{i\in \cN'} b'_i \hat{x}_i}}  b'_i \hat{x}_i\right\}.
\end{equation}

By \Cref{lem.sep}, $\cO$ is an exact polyhedral outer approximation, so replacing $x \in \cC$ with $x \in \cO$ does not change the binary feasible set. This gives rise to an exact MILP formulation equivalent to the submodular knapsack problem with conflicts:
    \begin{subequations}
    \label{pricing.ref.milp}
    \begin{align}
  		 \max \quad \sum_{i\in \cN'} \pi'_i x_i, \qquad& \\
   	 \st \qquad  \theta^\top x\le c, \qquad&  \forall \theta \in \Theta\\
   	 x_{i_1} + x_{i_2} \le 1, \qquad & \forall (i_1, i_2) \in \cD', \\
   	 x_i \in \{0,1\},\qquad & \forall  i  \in \cN'.
    \end{align}
    \end{subequations}

However, $\cX$ (and hence $\Theta$) is unknown before exploring the search space, and its cardinality may be exponential. In practice, the cuts corresponding to $\Theta$ can only be separated \textit{lazily}, i.e., a cut is added until a point $\hat{x}$ is found in $\cX$. Off-the-shelf solvers do not use this finite family of cuts, but it is a crucial component for constructing our PWL-B$\&$C algorithm in \Cref{sec.pwl}.

The following lemma explains the approximation error of the  polyhedral outer approximation $\cO$ w.r.t. $\overline{\cC}$.

\begin{lemma}[\cite{ben2001polyhedral}]
\label{lem.pout}
Let \(\epsilon > 0\), then there exists a method to construct a  polyhedral outer approximation $\cO$ of \( \overline{\cC}\) with additional \(\bO{1}|\cN'|\log(\frac{1}{\epsilon})\) variables and constraints, such that the relative \(\ell_\infty\) approximation error \(\max_{x \in \cO} |\sum_{i \in \cN'} a'_i x_i  + \sigma\sqrt{\sum_{i\in \cN'} b'_i x_i^2 } - c| / c\) is at most \(\epsilon\).
\end{lemma}

Note that the approximation error of the polyhedral outer approximation depends on the number of variables.

\subsection{MBQCP formulation}

We present a non-convex Mixed Binary Quadratically Constrained Programming formulation for the submodular knapsack problem (with conflicts). Although we do not use this formulation to solve the price subproblems, this formulation inspires PWL relaxation and the PWL-B$\&$C algorithm. Here, we introduce a slack variable $w$ to define the sum $\sum_{i \in \cN'} a'_i x_i$. Then our MBQCP formulation becomes the following non-convex MINLP program:
\begin{subequations}
\label{pricing.ref.mbqcp}
\begin{align}
  	 \max \quad \sum_{i\in \cN'} \pi'_i x_i, \qquad& \\
    \st \qquad  \sum_{i \in \cN'} a'_i x_i = w, \qquad & \\
    \sigma^2\sum_{i \in \cN'} b'_i x_i  \le (c - w)^2,\qquad& \label{pricing.ref.mbqcp.quad}\\
    x_{i_1} + x_{i_2} \le 1, \qquad & \forall (i_1, i_2) \in \cD', \\
    x_i \in \{0,1\},\qquad & \forall  i  \in \cN', \\
    w \in [0,c].\qquad &
\end{align}
\end{subequations}
Although the program contains a concave quadratic constraint \eqref{pricing.ref.mbqcp.quad}, the nonlinearity is only a univariate quadratic function compared to the \(|\cN'|\)-dimensional nonlinear SOC function $f$ in \eqref{pricing.ref.bsocp.cap}.

\subsection{PWL relaxation}
\label{sec.pwlrel}
A Piece-Wise Linear (PWL) function is linear on each piece of a given partition of its domain. We derive a MILP relaxation of the MBQCP formulation \eqref{pricing.ref.mbqcp} based on the PWL relaxation for the quadratic function, and refer to this new MILP relaxation as the PWL relaxation. The approximation error of the optimal PWL relaxation is discussed in this section. Let us denote by \(q(w):=(c-w)^2\) the univariate quadratic function. We denote a value of the slack variable \(w\) in the constraint \eqref{pricing.ref.mbqcp.quad} as a \textit{breakpoint}. Given an ordered set of breakpoints \(\cB = (w_1, w_2,\hdots, w_h)\) such that \(w_k \in [\underline{w},\overline{w}]\) ($k \in [h]:=\{1,\cdots,h\}$), $w_1 = \underline{w}$ and $w_h = \overline{w}$, the following function is a PWL approximation of \(q\) over the domain \([\underline{w},\overline{w}]\):
\begin{equation*}
   \bar{q}_{\cB}(w):=
   	 \frac{q(w_k)-q(w_{k-1})}{w_{k}-w_{k-1}}(w-w_{k-1}) + q(w_{k-1}), \textup{ for } w_{k-1}\leq w\leq w_{k}, 2 \le k \le h.
\end{equation*}
Note that $ \bar{q}_{\cB}$ is  an \textit{over-estimator} of  $q$ due to the convexity of $q$.

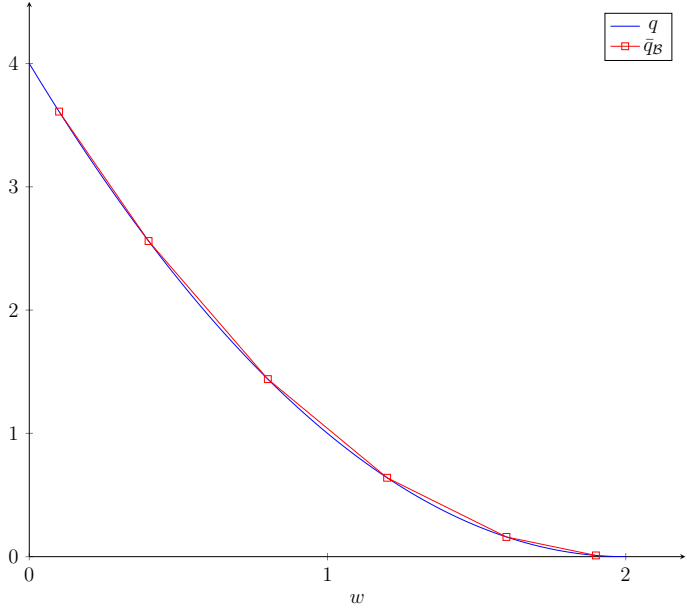
\begin{figure}
    \centering
\begin{tikzpicture}
[
thick,scale=0.65,
 declare function={
 quad(\x)=(2-\x)*(2-\x);
 secant(\a,\b,\x)=and(\x in \[\a,\b\])*(quad(\b)-quad(\a))/(\b-\a)*(\x-\a)+quad(\a);
 }
]
\begin{axis}[
axis lines = left,
xmax=2.2,xmin=0,
ymin=0,ymax=4.5,
xlabel=$w$,ylabel=,
xtick={0,...,2},
ytick={0,...,4},
width=15cm,
anchor=center,
]
\addplot [
    domain=0:2,
    samples=200,
    color=blue,
]{quad(x)} ;
\addplot[
color=red,
mark=square,
]
coordinates {
(0.1, {quad(0.1)})(0.4,{quad(0.4)})(0.8,{quad(0.8)})(1.2,{quad(1.2)})(1.6,{quad(1.6)})(1.9,{quad(1.9)})};
 \addlegendentry[]{$q$}
 \addlegendentry[]{$\bar{q}_{\cB}$}
\end{axis}
\end{tikzpicture}
    \caption{Graphs of the quadratic and its PWL over-estimator}
    \label{fig:func}
\end{figure}

We call \(\cB\) a breakpoint set in \([\underline{w},\overline{w}]\), and \(\bar{q}_{\cB}\) its  induced PWL function. Note that we consider the two bounds  \(\underline{w}\) and \(\overline{w}\) as breakpoints here. \Cref{fig:func} shows the graphs of a quadratic function and its PWL over-estimator, where $\underline{w}=0.1$, $\overline{w}=1.9$, $c=2$, and $\cB = \{0.1, 0.4,0.8,1.2,1.6, 1.9\}$.

Assume that we are given the breakpoints $\cB$. Replacing $  \sigma^2\sum_{i \in \cN'} b'_i x_i \le q(w)$ with  $  \sigma^2\sum_{i \in \cN'} b'_i x_i \le \bar{q}_\cB(w)$ in the constraint \eqref{pricing.ref.mbqcp.quad}, we obtain the following PWL relaxation of the MBQCP formulation \eqref{pricing.ref.mbqcp}:

\begin{subequations}
\label{pricing.rel.mbqcp}
\begin{align}
  	 \max \quad \sum_{i\in \cN'} \pi'_i x_i, \qquad& \\
    \st \qquad  \sum_{i \in \cN'} a'_i x_i = w, \qquad & \label{pricing.rel.mbqcp.eq} \\
    \sigma^2\sum_{i \in \cN'} b'_i x_i  \le  \bar{q}_{\cB}(w),\qquad& \label{pricing.rel.mbqcp.quad}\\
    x_{i_1} + x_{i_2} \le 1, \qquad & \forall (i_1, i_2) \in \cD', \\
    x_i \in \{0,1\},\qquad & \forall  i  \in \cN', \\
    w \in [0,c].\qquad &
\end{align}
\end{subequations}

\begin{remark}
(Modeling PWL functions) The graphs of PWL functions have several MILP formulations, see \cite{vielma2010}. In this paper, we consider the logarithmic model. We denote by $z$ the auxiliary binary variables introduced in the MILP formulation of $\bar{q}_{\cB}$. From version 20.1.0 \citep{cplexsos}, \texttt{CPLEX} can automatically formulate $\bar{q}_{\cB}$ to the logarithmic model and add auxiliary variables $z$ in the internal data structure.
\end{remark}

The approximation error of a PWL relaxation is expressed as $\ell_p$-norm of the difference between the approximation function and the target function.

\begin{definition}
Given a set \(\cB \subset [\underline{w}, \overline{w}]\) of breakpoints, the \(\ell_p\) approximation error of \(\bar{q}_{\cB}\) with respect to \(q\) over \([\underline{w},\overline{w}]\) is defined as \(\ell_p(\bar{q}_{\cB}, q) := (\int_{ \underline{w}}^{\overline{w}}|\bar{q}_{\cB}(w) - q(w)|^p \,dw)^{\frac{1}{p}}. \)
\end{definition}

Since the approximation error measures the quality of a PWL approximation to the quadratic function, thus it in turn measures the error of the PWL relaxation to the MBQCP formulation. Empirically, the optimal solution to a PWL relaxation with a small approximation error should have a small gap to the optimal solution of the submodular knapsack with conflicts.  On the other hand, although adding breakpoints decreases the approximation error, it increases the computation resource to solve the PWL relaxation. So a common problem is understanding the best possible approximation error given a fixed number of breakpoints (limited computational resource).

This problem can be formalized as follows. Given an integer \(h\) (number of breakpoints), denote by \(\bB^h\) the family of breakpoint sets of cardinality \(h\) in \([\underline{w}, \overline{w}]\), the \textit{breakpoint selection problem} aims to find a set \(\cB \in \bB^h\) to minimize the  \(\ell_p\) error:

\begin{equation}
\label{bsp}
    \min_{\cB \in \bB^h} \ell_p(\bar{q}_{\cB}, q).
\end{equation}

 A convex program
\cite{geissler2012using} can compute the $\ell_\infty$-approximation error for general nonconvex functions.  An error analysis \cite{berjon2015optimal} gives asymptotically tight bounds to quantify the $\ell_2$-approximation error.

The following theorem gives the best \(\ell_\infty\)-approximation error that we can achieve: an optimal solution to the breakpoint selection problem under the $\ell_\infty$-approximation error is an equidistant partition of \([\underline{w}, \overline{w}]\).

\begin{theorem}
\label{lem.approx}
Given $\cB \in \bB^h$,
 $$\ell_\infty(\bar{q}_{\cB}, q) = \max_{w \in [\underline{w},\overline{w}] }|\bar{q}_{\cB}(w) - q(w)|= \max_{2\le k \le h} \frac{(w_k - w_{k-1})^2}{4}.$$ Furthermore, let \(w_k = \underline{w} + \frac{k-1}{h-1}(  \overline{w} - \underline{w})\) for \(1 \le k \le h\), which yields the minimum \(\ell_\infty\)-approximation error \(\frac{(\overline{w} - \underline{w})^2}{4(h-1)^2}\) for the break point selection problem \eqref{bsp}.
\end{theorem}

The  proof can be found in \ref{append.lem.approx}.
The approximation error decreases with the quadratic rate with respect to \(h\). The relative \(\ell_\infty\)-approximation error is defined as $$\frac{\ell_\infty(\bar{q}_{\cB}, q)}{(  \overline{w} - \underline{w})^2}.$$

We have the following result on the relative approximation error of the PWL relaxation.
\begin{corollary}
\label{cor.pout}
Let \(\epsilon > 0\), then there exists a MILP formulation of PWL function $\bar{q}_{\cB}$ induced by $\cB$ with \(\bO{1}\log(\frac{1}{\epsilon})\)  binary variables and \(\bO{1}\frac{1}{\sqrt{\epsilon}}\)  continuous variables and constraints, such that the relative $\ell_\infty$-approximation error is at most \( \epsilon\).
\end{corollary}
\begin{proof}
For the logarithmic model of PWL function, given $h$ breakpoints from the equidistant partition, the relative \(\ell_\infty\)-approximation error is $\frac{(\overline{w} - \underline{w})^2}{4(h-1)^2(\overline{w} - \underline{w})^2}=\frac{1}{4(h-1)^2}$ with $\log(h-1)$ binary variables and $h-1$ continuous variables and constraints \citep{vielma2010}, the result follows.
 \end{proof}

Next, we summarize the approximation errors of two relaxations to their corresponding formulations. Note that we do not consider the integrality of the binary variable $x'$. Comparing \Cref{lem.pout} and \Cref{cor.pout}, the approximation error of the PWL relaxation \eqref{pricing.rel.mbqcp} to the MBQCP formulation \eqref{pricing.ref.mbqcp} is independent of the number of variables, while the approximation error of the polyhedral outer approximation to the BSOCP formulation \eqref{pricing.ref.bsocp} depends on this number.

We remark that our PWL relaxation differs from \cite{shylo2013}'s PWL relaxation. The constraint \eqref{pricing.ref.mbqcp} of MBQCP formulation is equivalent to $\sigma \sqrt{\sum_{i \in \cN'} b'_i x_i}  \le c -w$, and a  PWL relaxation was used for  the left-hand side concave function $\sigma \sqrt{\sum_{i \in \cN'} b'_i x_i}$  in \cite{shylo2013}. However, the optimal approximation error for such PWL relaxation has yet to be discovered.

\subsection{Exact PWL-B\&C algorithm}
\label{sec.pwl}
The approximation error of the PWL relaxation is dimensionless but only for a small number of breakpoints, it is not exact. Instead of adding many breakpoints, the finite number of cuts induced by the set $\Theta$ in \eqref{corecuts} suffices to make the PWL relaxation exact. To solve it, we propose a combined formulation and a branch-and-cut algorithm based on the PWL relaxation (PWL-B\&C).

\begin{subequations}
\label{pricing.rel.comb}
\begin{align}
  	 \max \quad \sum_{i\in \cN'} \pi'_i x_i, \qquad& \\
    \st \qquad  \sum_{i \in \cN'} a'_i x_i = w, \qquad & \label{pricing.rel.comb.eq} \\
    \sigma^2\sum_{i \in \cN'} b'_i x_i  \le  \bar{q}_{\cB}(w),\qquad& \label{pricing.rel.comb.quad}\\
    \qquad \theta ^\top  x\le c,\qquad & \forall \theta \in \Theta \label{pricing.ref.comb.cap}\\
    x_{i_1} + x_{i_2} \le 1, \qquad & \forall (i_1, i_2) \in \cD', \\
    x_i \in \{0,1\},\qquad & \forall  i  \in \cN', \\
    w \in [0,c].\qquad &
\end{align}
\end{subequations}
Formulation \eqref{pricing.rel.comb} combines the MILP formulation \eqref{pricing.ref.milp} with the (redundant) PWL relaxation. As already mentioned by \Cref{lem.sep}, the MILP formulation \eqref{pricing.ref.milp} is an exact formulation for submodular knapsack problems with conflicts, so this combined formulation \eqref{pricing.rel.comb} is also exact.

The intuition underlying the combined formulation \eqref{pricing.rel.comb} is that we cannot add numerous valid inequalities \eqref{pricing.ref.comb.cap} \textit{a priori}. In practice, we add them \textit{lazily} to exclude infeasible binary solutions to the submodular knapsack with conflicts in the course of the search, and this method is typically supported or suggested by \textit{lazy cut callbacks} of some solvers such as \texttt{CPLEX} and \texttt{SCIP}. However, in this way,  we cannot control the initial relaxation quality given solely by a few  valid inequalities from \eqref{pricing.ref.comb.cap}. On the contrary, the PWL relaxation \eqref{pricing.rel.comb.quad} can be enforced \textit{a priori}, and its quality is controllable (\Cref{lem.approx}). So we can leverage it to reduce the initial search space and refine the relaxation by adding valid inequalities lazily. This intuition and formulation give rise to a tailored \Cref{algo} for submodular knapsack with conflicts, partly inspired by algorithms in \cite{coey2018outer}. We show in experiments that this formulation with redundant constraints \eqref{pricing.rel.comb.eq} and \eqref{pricing.rel.comb.quad} can be solved much faster than the standard BSOCP formulation \eqref{pricing.ref.bsocp}. In practice, only a few cuts in \eqref{pricing.ref.comb.cap} must separate before the convergence.

Our algorithm consists of three main steps: tightening the bounds, constructing the PWL relaxation (breakpoints), and the PWL B\&C algorithm. First, bound tightening is a preprocessing procedure used to tighten the bounds on the breakpoints for all pricing problems. Then, the PWL relaxation (breakpoints) is constructed for all pricing problems, and this is also a pre-solving procedure. The construction depends on the number of items, the size of the items, and the capacity. Finally, based on the PWL relaxation, the PWL-B\&C algorithm is adapted to the LP-B\&C algorithm \citep{coey2018outer}.\\

\textbf{Bound tightening} The bound tightening procedure is called before the branch-and-price algorithm to shrink the boundaries of the breakpoints \(\cB\) into \([0,c]\).
 
Considering a pricing problem at a node of the search tree, we find that if $w = \sum_{i \in \cN'} a'_i x_i$ is small, $q(w)=(c-w)^2$ is larger than $\sigma^2 \sum_{i\in \cN'} b'_i x_i$, so the capacity constraint \eqref{pricing.ref.mbqcp.quad} is not active. Thus, there is no need to overestimate $q$ when $w$ is small. More precisely, there is a \(\underline{w} \in [0,c]\) such that, for any binary solution \(x \in \{0,1\}^{\cN'}\), let \( w = \sum_{i \in \cN'} a'_i x_i \), if \(w \le \underline{w}\), then \(\sigma^2\sum_{i \in \cN'} b'_i x_i \le q(\underline{w})\). Since $q$ is non-increasing, \(q(w) \ge q(\underline{w}) \ge \sigma^2\sum_{i \in \cN'} b'_i x_i\). The point $\underline{w}$ is called \textit{lower breakpoint}, the submodular capacity constraint \eqref{pricing.ref.mbqcp.quad} is never violated for $w \in [0,\underline{w}]$. We can start by overestimating \(q\) starting from the maximum lower breakpoint computed from the following convex MBQCP problem:
 
 \begin{equation}
\label{bdt.y1}
\begin{split}
	 \underline{w} := \max \quad w, \qquad& \\
    \st \qquad  \sum_{i \in \cN'} a'_i x_i = w,\qquad& \\   
   \qquad \sigma^2 \sum_{i\in \cN'} b'_i x_i \ge (c - w)^2,\qquad& \\
    x_{i_1} + x_{i_2} \le 1, \qquad & \forall (i_1, i_2) \in \cD', \\
    x_i \in \{0,1\},\qquad & \forall  i  \in \cN'.
\end{split}
\end{equation}

Similarly, we can define the \textit{upper breakpoint}.
There exists some upper breakpoint \(\overline{w} \in [0,c]\), such that, for every binary solution \(x \in \{0,1\}^{\cN'}\), if \(\sum_{i \in \cN'} a'_i x_i  + \sigma\sqrt{\sum_{i\in \cN'} b'_i x_i^2 } \le c\), then \( \sum_{i \in \cN'} a'_i x_i  \le \overline{w} \). The minimum upper breakpoint can be computed from the following BSOCP problem:

 \begin{equation}
 \label{bdt.ym}
\begin{split}
   	 \overline{w}:= \max \quad \sum_{i\in \cN'}a'_i x_i , \qquad& \\
    \st \qquad  \sum_{i \in \cN'} a'_i x_i  + \sigma\sqrt{\sum_{i\in \cN'} b'_i x_i^2 } \le c,\qquad& \\
    x_{i_1} + x_{i_2} \le 1, \qquad & \forall (i_1, i_2) \in \cD', \\
    x_i \in \{0,1\},\qquad & \forall  i  \in \cN'.
\end{split}
\end{equation}

We solve the above two programs at the root node and obtain the bound \([\underline{w}_r, \overline{w}_r]\) for breakpoints. Since the feasible sets of the other nodes are a subset of the root node set, the above programs at other nodes are more strict than those at the root node. It follows for $\underline{w}, \overline{w}$ of any other node that \(\underline{w}_r \le \underline{w}\) and \( \overline{w} \le \overline{w}_r\). We then set $\underline{w} = \underline{w}_r$ and $\overline{w} = \overline{w}_r$ for all nodes.\\

\textbf{Construction of breakpoints} To determine the number of breakpoints \(\cB\), we run a greedy heuristic algorithm that tries to maximize the number of items in a bin. We take $h$ as the solution value given by the heuristic algorithm and assign breakpoints \(h\) equidistantly in \([\underline{w}, \overline{w}]\). The equidistant partition gives the best approximation error according to \Cref{lem.approx} for a fixed number of breakpoints. We also add a breakpoint corresponding to $w=0$.\\

 \textbf{PWL-B\&C algorithm} The main steps of the PWL-B\&C algorithm are described in \Cref{algo}. Recall that the problem \eqref{pricing.ref} is a maximization problem. \Cref{algo} maintains  a set of active nodes $\mathscr{N}$ of the search tree, a pool of cuts \(\mathscr{C}\),  an incumbent solution $x^\ast$ ($\sum_{i \in \cN'}\pi'_ix^\ast_i$ is a primal bound).
 
A node $(l, u, U)$ is characterized by the finite variable boundary vectors $l$ and $u$ and the node's dual upper bound $U$. The upper bound $U$ is inherited from its parent node and computed via the LP relaxation. Note that the PWL function is a modeling concept. We use a MILP solver, i.e., \texttt{CPLEX}, that formulates the PWL function $q_{\cB}$ into a MILP. We denote by $z$ the additional binary variables to model \(\bar{q}_{\cB}\) (see \Cref{sec.pwlrel}).  The variables $z$ are also constructed internally by \texttt{CPLEX}, and we assume that the PWL function is forced when $z$ is set to binary.

We denote by \(\mathscr{M}_\cB(\mathscr{C}, l, u, U)\) the MILP relaxation restricted to finite bounds $(l,u)$ for $(x,z)$ at a node of the search tree. The MILP relaxation \(\mathscr{M}_\cB(\mathscr{C}, l, u, U)\) consists of the PWL relaxation \eqref{pricing.rel.mbqcp}, cuts from $\mathscr{C}$, and other cuts added by the MILP solver.

\begin{algorithm}[htbp]
\SetAlgoLined
 \textbf{Input:} a submodular knapsack problem with conflicts \eqref{pricing.rel.comb}, and the set $\cB$ of breakpoints\;
   \textbf{Output:} a primal solution $x^\ast$ and a dual upper bound (dual gap)\;
    initialize MILP \(\mathscr{M}_\cB(\mathscr{C}, l_0, u_0, \infty)\) as the PWL relaxation \eqref{pricing.rel.mbqcp} \Comment*[r]{ the PWL function is modeled by auxiliary binary variable $z$}
  initialize  cut pool \(\mathscr{C}\) to \(\emptyset\), the node list \(\mathscr{N}\) of  \(\mathscr{M}_\cB(\mathscr{C}, l_0, u_0, \infty)\) with root node \((l_0,u_0)\), incumbent solution \(x^\ast = 0\), and the  upper bound  \(U\)  of the root node to \(\infty\)\;
 \While{\(\mathscr{N}\) contains nodes}{
  remove a node \((l,u)\) from \(\mathscr{N}\) \; \label{algo:rmv}
   solve LP relaxation of \(\mathscr{M}_\cB(\mathscr{C},l,u, U)\)\; \label{algo:sollp}
   \uIf{ LP is infeasible }{
   \Continue
   \Comment*[r]{fathomed by infeasibility} \label{algo:lpinf}
   }
   get an LP optimal solution \((\hat{x}, \hat{z})\)\;
   \label{algo:calcu}
   \uIf{ upper bound \(U \le \sum_{i \in \cN'} \pi_i \hat{x}_i\)}{
   \Continue \Comment*[r]{fathomed by bound} \label{algo:fbd}
   }
   \Else{
   set $U$ to $\sum_{i \in \cN'} \pi_i \hat{x}_i$ \Comment*[r]{update the dual upper bound} \label{algo:upbd}
   }
   \If{ \((\hat{x},\hat{z})\) is binary}{
  	 \uIf{  \(\hat{x}\) satisfies capacity constraint \eqref{pricing.ref.cap}}{
 		 set \(x^\ast\) to  \(\hat{x}\)\; \label{algo:upsol}
   		 \Continue  \Comment*[r]{fathomed by integrality} \label{algo:fint}
  	 }
   	 \Else{
   		 add separation cut to \(\mathscr{C}\) by \Cref{lem.sep}\; \label{algo:sep}
   		 add the node \(\mathscr{M}_\cB(\mathscr{C},l,u, U)\) to \(\mathscr{N}\) \; \label{algo:cut}
   		 \Continue  \Comment*[r]{reoptimization after cut added} \label{algo:reopt}
   	 }
   }
   add branch nodes to \(\mathscr{N}\) using \((\hat{x}, \hat{z})\) (fractional) and \(U\)\; \label{algo:brnch}
 }
 \caption{PWL-B$\&$C algorithm}
 \label{algo}
\end{algorithm}

The node set $\mathscr{N}$ initially contains the root node $(l^0, u^0)$, where $l^0, u^0 \in \bR^I$ are the finite initial global bounds on variables \((x,z)\).
On \Cref{algo:rmv} of \Cref{algo}, the main loop removes a node $(l, u, U)$ from $\mathscr{N}$.  \Cref{algo:sollp} solves the  LP  relaxation  of \(\mathscr{M}_\cB(\mathscr{C}, l, u, U)\) by the MILP solver.

If the LP-relaxation \(\mathscr{M}_\cB(\mathscr{C}, l, u, U)\) is infeasible, \Cref{algo:lpinf} immediately fathoms the node by infeasibility. The upper bound $U$ of the node means that any feasible solution to the combined formulation \eqref{pricing.rel.comb} that satisfies the bounds of the node for \((x,z)\) has an objective value of at most $U$. Since LP is a relaxation of the combined formulation \eqref{pricing.rel.comb}, any feasible solution to the combined formulation \eqref{pricing.rel.comb} that satisfies the bounds of the node for \(x\) has a objective value of at most $U$.

\Cref{algo:fbd} fathoms the node by bound if $U$ is not better than the incumbent value. Otherwise, the upper bound $U$ of the node is set to the optimal value of LP on \Cref{algo:upbd}.

If $\hat{z}$ is not binary, then the PWL function is not implicitly enforced by the integrality of $\hat{z}$, so the algorithm should continue to branch. If $\hat{z}$ is binary (the PWL function is enforced) and $\hat{x}$ is binary, then the algorithm examines the solution $\hat{x}$.

If additionally, $\hat{x}$ is feasible (the capacity constraint is satisfied), then its objective value  should be at least the upper bound \(U\). \Cref{algo:upsol} stores the new incumbent solution $\hat{x}$, and \Cref{algo:fint} fathoms the node since \(\hat{x}\) is an optimal binary solution with respect to the bounds $(l,u)$.  Otherwise, the constraint \eqref{pricing.ref.comb.cap} is  violated. \Cref{algo:sep} adds this constraint to the cut pool \(\mathscr{C}\), \Cref{algo:cut} adds the current node for re-optimization, and \Cref{algo:reopt} discards $\hat{x}$ by the cut in the next optimization iteration. Finally, $(\hat{x},\hat{z})$ must be fractional on \Cref{algo:brnch}, the algorithm branches using the information from fractionality and $U$.

We remark that the idea in \cite{shylo2013} for exact algorithms does not deploy  cutting planes for approximating submodular knapsack, so in each iteration, new breakpoints are added to PWL, relaxations, and the underlying MILP solver needs restarts.

\subsection{Hybrid pricing strategy}
The pricing heuristic in \Cref{sec.heur} is fast, but it cannot guarantee the dual upper bound required by the Farley bound of \Cref{lem.valildlbd}. The exact pricing algorithm is slow but yields the dual upper bound for the pricing problem. The hybrid pricing strategy first calls the pricing heuristic  to decide whether the exact pricing algorithm can improve the local dual bound of the master problem.

In fact, the exact algorithm is required only under a particular condition. The following proposition gives the condition.
 
\begin{proposition}
\label{lem.bdimprove}
Let \(v_{heur}\) be the solution value of the pricing heuristic, let \(v_{\mathrm{RMLP}}\) be the optimum of RMLP \eqref{master.dw.rel}, and let \(v_{\mathrm{ld}}\) be the current local dual bound for the master problem. If \(\frac{v_{\mathrm{RMLP}}}{v_{heur}} \le v_{\mathrm{ld}}\), the exact algorithm cannot yield a better local dual bound than $v_{\mathrm{ld}}$.
\end{proposition}
\begin{proof}
Let \(v_{popt}\) be the optimum for the pricing problem \eqref{pricing.ref}, then \(v_{heur} \le v_{popt}\). It follows that \(\frac{v_{\mathrm{RMLP}}}{v_{popt}}\le \frac{v_{\mathrm{RMLP}}}{v_{heur}} \le  v_{\mathrm{ld}}\). However, $v_{popt}$ is the smallest pricing dual bound $v_{\mathrm{price}}$, so $\frac{v_{\mathrm{RMLP}}}{v_{popt}}$ is the greatest Farley bound according to \Cref{lem.valildlbd}. Therefore even if the pricing algorithm is solved to optimality, we cannot obtain a better bound than \(v_{\mathrm{ld}}\).
 \end{proof}

 If the condition \(\frac{v_{\mathrm{RMLP}}}{v_{heur}} \le v_{\mathrm{ld}}\) holds, one can get rid of the exact pricing algorithm, and use the solution from the fixing-greedy heuristic in \Cref{sec.heur}. The hybrid pricing strategy is outlined in Algorithm \ref{sched}.

\begin{algorithm}[htbp]
\SetAlgoLined
 \textbf{Input:} a pricing problem  \eqref{pricing.ref} with the objective coefficients $\pi'$,  \(v_{\mathrm{RMLP}}\) the optimum of RMLP \eqref{master.dw.rel},  \(v_{\mathrm{ld}}\) the local dual bound of the master problem\;
   \textbf{Output:} a generated column $x^\ast$,  and the updated local dual bound $v_{\mathrm{ld}}$\;
  \label{sched:updateweights}
 call the pricing heuristic with the objective coefficients  $ \pi'$ \Comment*[r]{run heuristic first} \label{sched:runheur1}
 let $x, v_{heur}$  be the heuristic solution and its value\;
 \uIf{$\frac{v_{\mathrm{RMLP}}}{v_{heur}} \le v_{\mathrm{ld}}$ and $1 - \sum_{i\in \cN'} \pi'_i \bar{x}_i < 0$}  
 {

    $x^\ast \leftarrow x$ \Comment*[r]{heuristic solution} \label{sched:originsol}
 }
 \Else{
   call the exact pricing Algorithm \ref{algo}  \Comment*[r]{exact pricing} \label{sched:exact}
  let $\tilde{x}, v_{\mathrm{price}}$  be the primal solution and the dual bound\;
  $x^\ast \leftarrow \tilde{x}$\;
   $v_{\mathrm{ld}} = \max \{v_{\mathrm{ld}}, \frac{v_{\mathrm{RMLP}}}{v_{\mathrm{price}}}\}$   \Comment*[r]{update the local dual bound} \label{sched:updatedb}
 }
 \caption{Hybrid pricing strategy}
 \label{sched}
\end{algorithm}

 The heuristic algorithm is called first in Line \ref{sched:runheur1}. If  $\frac{v_{\mathrm{RMLP}}}{v_{heur}} \le v_{\mathrm{ld}}$, the exact pricing is not needed.  If the heuristic solution $x$ has a negative reduced cost, the strategy outputs it in Line \ref{sched:originsol}. Otherwise, the strategy calls the exact algorithm in Line \ref{sched:exact}.

\section{Computational experiments}
\label{sec.experi}
In this section, we present the computational experiments we made to test the effectiveness of our branch-and-price algorithms for SMBP. In particular, we test different configurations of branch-and-price algorithms to evaluate the proposed techniques. The source code and benchmarks are publicly available on the project website \href{https://github.com/lidingxu/cbp}{https://github.com/lidingxu/cbp}. We also provide a bash file to reproduce the experiments on Linux systems.
\subsection{Benchmarks}
We produce benchmarks as described in \cite{cohen2019overcommitment}. The authors test their approximation algorithms on benchmarks from real cloud data centers of \texttt{Google}, which are not accessible due to confidentiality \footnote{We, therefore, create new instances using the same generation method.}.

They also describe data generation methods by considering a variety of uncertainty models, and these methods have a probabilistic interpretation: parameters of SMBP instances are derived from parameters of uncertainty models. For different risk levels $\alpha$, they propose three data generation methods (cases) to construct the data $a,b, \sigma$ in SMBP \eqref{master.submod}, i.e., the Gaussian case, the Hoeffding inequality case, and the distributionally robust approximation case.

 We describe the generation methods next. In summary, we first determine overall parameters such as capacity and item numbers, then we generate distributions, and finally cast  parameters of distributions into parameters of items.

The overall parameters of instances are set as follows.
We set the capacity of each bin to 72 (the number of cores of the servers),  the risk level $\alpha \in \{0.6, 0.7, 0.8, 0.9, 0.95,0.99\}$. We set the number of items (i.e., jobs) $|\cN| \in \{100, 400, 1000\}$ to obtain three benchmarks with different sizes: \texttt{CloudSmall}, \texttt{CloudMedium}, and \texttt{CloudLarge}. There are three generation methods and six risk levels.

The distributions of instances are set as follows.
We call the distribution of $\mu_i$ the \textit{target distribution} for item $i$. We assume that every $\mu_i$ follows the same target distribution. This target distribution is  unknown in \cite{cohen2019overcommitment} except for  its quantiles in \Cref{itemdist}.

Given $\alpha$ and $\cN$, we generate an SMBP instance as follows:
\begin{enumerate}
    \item sample $\mu_i$  ($i \in \cN$)  according to  \Cref{itemdist};
    \item sample $a$ and $b$ from $\mu$ and $\sigma$, using one of the following cases:
    \begin{itemize}
    \item Gaussian case;
    \item Hoeffding’s inequality case;
    \item distributionally robust approximation case.
    \end{itemize}
\end{enumerate}

\begin{table}[htbp]
\centering
\caption{Example distribution of item size}
\begin{tabular}{c c c c c c c c}
\hline
 Item sizes &  1 &  2 &  4 & 8 & 16 & 32  & 72\\
 \hline
\% Items  & 36.3  & 13.8 & 21.3 & 23.1 & 3.5 & 1.9 & 0.1\\
\hline
\end{tabular}
\label{itemdist}
\end{table}

We first illustrate the approach of sampling $\mu$. We approximate the target distribution by a normalized histogram such that its quantile distribution is the same as in \Cref{itemdist}. A histogram consists of intervals divided from the entire range $[0,72]$, and each interval has endpoints of two consecutive quantiles of \Cref{itemdist}. The histogram gives a discrete non-parametric estimation of the target distribution. We apply a two-stage sampling to obtain a nominal item size $\mu_i$ ($ i\in \cN$) sampled from a continuous distribution. It has two steps:

\begin{enumerate}
    \item sample an interval $[d_1,d_2]$ from the histogram;
    \item sample a nominal item size $\mu_i$ from $[d_1,d_2]$ uniformly.
\end{enumerate}

Second, we construct a truncated Gaussian, which is defined by its lower and upper bounds $\underline{A}$ and $\overline{A}$, mean  $\mu'$, and its standard deviation $\sigma'$. To obtain these parameters, for each $i \in \cN$,  we:

\begin{enumerate}
    \item sample $\underline{A}_i \in [0.3, 0.6]$ and $\overline{A}_i \in [0.7, 1.0]$ uniformly;
    \item sample scale parameter $s_i \in [0.1, 0.5]$;
    \item compute the mean $\mu'_i$ and the standard variation $\sigma'_i$ of the truncated Gaussian with lower bound $\underline{A}_i$, upper bound $\overline{A}_i$ and scale parameter $s_i$.
\end{enumerate}

With the above parameters, we generate the data  $a,b, \sigma$ of SMBP \eqref{master.submod}. There are three cases, which correspond to different assumptions on the uncertainty or probability distribution.

For the Gaussian case:
\begin{enumerate}
    \item let $\sigma = \Phi^{-1}(\alpha)$, where $\Phi$ is the cumulative distribution function  of the Gaussian distribution;
    \item for $i \in \cN$, let $a_i = \mu'_i \mu_i$ and $b_i = {(\sigma'_i \mu_i)}^2$.
\end{enumerate}

For the Hoeffding’s inequality case:
\begin{enumerate}
    \item let $\sigma = \sqrt{-0.5\ln{(1- \alpha)}}$;
    \item for $i \in \cN$, let $a_i = \mu'_i \mu_i$ and $b_i = {((\overline{A}_i-\underline{A}_i) \mu_i)}^2$
\end{enumerate}

For the distributionally robust approximation case:
\begin{enumerate}
    \item let $\sigma = \sqrt{\alpha/(1-\alpha)}$;
    \item for $i \in \cN$, let $a_i = \mu'_i \mu_i$ and $b_i = {(\sigma'_i \mu_i)}^2$.
\end{enumerate}

For all the above cases, if there exists $i \in \cN$ such that $a_i, b_i$ are too large to fit a bin (usually for large $\alpha, \sigma$), then we rescale $a_i, b_i$ to fit the bin.

 We generate six instances with different random seeds for each combination of generation methods and risk levels. As a result, we have $108 = 6 \times 6 \times 3$ instances in a benchmark.

\subsection{Experimental setups}
In this section, we describe the setup of the experiments, including the development environment, the implementation of the algorithms, and the solution statistics.\\

\textbf{Development environment} The experiments are conducted on a server with  Intel Xeon W-2245 CPU @ 3.90GHz, 126GB main memory, and Ubuntu 18.04 system. We use \texttt{SCIP} 8.0.1 \citep{gleixneretal2018oo} as a branch-and-price (B$\&$P) framework to solve the set cover formulation \eqref{master.dw}. We use ILOG \texttt{CPLEX} 22.1 as:
\begin{itemize}
    \item an LP solver to solve the RMLP \eqref{master.dw.rel};
    \item a BSOCP solver to solve the BSOCP formulations of SMBP  \eqref{master.submod} and the submodular knapsack problem with conflicts \eqref{pricing.ref.bsocp};
    \item a MILP solver used by the PWL-B$\&$C \Cref{algo};
\end{itemize}
\texttt{CPLEX}'s parameters are set by default, except we disable its parallelism. \\

\textbf{Solver implementation} We implement four solvers for SMBP according to the proposed techniques in this paper. Four of them are branch-and-price solvers. These solvers are as follows:
\begin{enumerate}
   \item  \bsocpcomp: a  solver using \texttt{CPLEX}'s B$\&$C algorithm to solve the compact BSOCP formulation of SMBP.
   \item  \dwbc: a B$\&$P solver for solving the set cover formulation \eqref{master.dw}, which uses  \texttt{CPLEX}'s B$\&$C algorithm to solve the BSOCP formulation \eqref{pricing.ref.bsocp} of the pricing problem.
   \item \dwpwl:  a B$\&$P solver for solving the set cover formulation \eqref{master.dw}, which uses   the PWL-B$\&$C algorithm to solve the combined formulation \eqref{pricing.rel.mbqcp} of the pricing problem.
   \item \dwhybrid: \dwpwl enhanced with the hybrid pricing strategy in \Cref{sched}.
\end{enumerate}

 We use the approximation algorithm from \cite{cohen2019overcommitment} to find an initial feasible solution that serves as a warm start for all solvers. All B$\&$P solvers deploy the column selection heuristic in \Cref{sec.pheur}, and all exact pricing solvers add the solution exclusion constraint \eqref{cons.exclude} to pricing problems. The time limit for each solver is 3600 CPU seconds.

If the column generation procedure at the root node does not finish after 3500 CPU seconds, it is halted, giving \texttt{SCIP} 100 CPU seconds to invoke its own primal heuristic.

For the pricing problems, we set the same time limit for the exact algorithms ($|\cN| \times 0.015$ CPU seconds) and the same tolerance for relative gaps. \\

\textbf{Performance metrics and statistical tests}  In order to evaluate the solver performance in
different instances, we compute  shifted geometric means (SGMs) (see \cite{achterberg2008constraint}) of performance metrics as aggregated statistics. Compared to arithmetic means, SMGs avoid the over-representation of biased outlier points. The SGM of values $v_1,...,v_N \geq 0$ with shift $s \geq 0$ is defined
as
\begin{equation*}
  \left(\prod_{i=1}^N (v_i + s)\right)^{1/N} - s.
\end{equation*}

Given an SMBP problem instance, let $\underline{v}$ be a dual lower bound
and $\overline{v}$ be a primal upper bound found by a solver. The relative dual gap in percentage is defined as:
\begin{equation*}
   \delta_d := \frac{\overline{v} - \underline{v}}{\overline{v}}\times 100.
\end{equation*}
A smaller relative dual gap indicates better performance.
 
 Let $v^a$ be the value of the solution found by the greedy min-utilization algorithm, which is communicated to all solvers as a warm start. The closed primal bound is defined as:
 \begin{equation*} \delta_p := \frac{v^a - \overline{v}}{ \max(\overline{v}- \underline{v}^\ast, 1e^{-6})} \times 100,
\end{equation*}
where $\underline{v}^\ast$ is the largest dual bound found among all solvers. A larger closed primal gap means better performance.

We report the following performance metrics for each instance tested by each solver and compute the SGMs of the benchmarks:
 
\begin{enumerate}
    \item $t$: the total running time in CPU seconds, with a shifted value set to 1;
   \item $\delta_d\%$: the relative dual gap in percentage, with a shifted value set to $1\%$;
   \item $\delta_p\%$: the closed primal bound in percentage, with a shifted value set to $1\%$;
   \item \#N:  the number of nodes of the search tree, with a shifted value set to 1;
  \item  \#C: the number of columns generated, with a shifted value set to $1$;
   \item E\%: the percentage of columns generated by the exact pricing algorithm, with a shifted value set to $1\%$;
   \item $\tau\%$: the relative dual gap in the percentage of a pricing problem solved by an exact algorithm, with a shifted value set to $1\%$;
	 \item $t_p\%$: the ratio between pricing time and total solving time in percentage, with a shifted value set to $1\%$.
\end{enumerate}

Metrics (1)-(4) refer to master problems and are available to all solvers. Metrics (5)-(8) refer to pricing problems and are not available for the \bsocpcomp, while metric (6) is $100\%$ for the \dwpwl and \dwbc.

\subsection{Comparative analysis of results}
\label{sec.mainresult}

The main computational results are summarized in \Cref{tab.ag}, for detailed results, we refer to \ref{append.detail}. For each benchmark, we report the SGM statistics of the performance metrics, the number of instances solved (denoted by \#S), and the number of instances with improved primal bounds (denoted by \#I). We also report  a computational test of adaptive selection of break points in \ref{sec.effect}. Next, we analyze the main computational results by comparing the solvers.

\begin{table}[]
\centering
\resizebox{0.99\columnwidth}{!}{
\begin{tabular}{|ll|*{6}{r}|*{4}{r}|}
\hline
\multicolumn{1}{|l}{\multirow{2}{*}{Benchmarks}} & \multicolumn{1}{|l|}{\multirow{2}{*}{Solvers}} & \multicolumn{6}{l|}{Problem statistics} & \multicolumn{4}{l|}{Pricing statistics} \\  \cline{3-9} \cline{10-12}
\multicolumn{1}{|l}{} & \multicolumn{1}{|l|}{}              		 &    $t$ & $\delta_d\%$ & $\delta_p\%$ & {\#N} & {\#S} & {\#I} &    {\#C} & E\% & {$\tau\%$} &$t_p \%$  \\
\hline
\multicolumn{1}{|l|}{\multirow{5}{*}{\shortstack{\texttt{CloudSmall}\\ ($|\cN|=100$)}}} & \bsocpcomp  &  1452 & 15.8 & 0.0 & 26601 & 18 & 0 &  - & - & - & -    \\
\multicolumn{1}{|l|}{}  & \dwbc &   2129 & 11.4 & 0.9 & 21 & 20 & 17 &    1373 & 100 & 3.56 & 99  	 \\
\multicolumn{1}{|l|}{} & \dwpwl &    633 & 2.4 & 2.7 & 66 & 61 & 32 &     1869 & 100 & 0.01 & 99  	 \\
\multicolumn{1}{|l|}{} & \dwhybrid &  330 & 2.0 & 3.4 & 127 & 65 & 36 &  	 3485 & 18 & 0.01 & 96 		 \\
\hline
\multicolumn{1}{|l|}{\multirow{5}{*}{\shortstack{\texttt{CloudMedium}\\ ($|\cN|=400$)}}} & \bsocpcomp  &    3600 & 100.0 & 0.0 & 0 & 0 & 0 & 	 -   &   - 	 &  -  	 &    -     \\
\multicolumn{1}{|l|}{}& \dwbc &    3600 & 39.0 & 0.1 & 2 & 0 & 4 &   861 & 100 & 0.39 & 98   	 \\
\multicolumn{1}{|l|}{}& \dwpwl &    3600 & 17.2 & 0.4 & 1 & 0 & 10 &   3372 & 100 & 0.01 & 91     \\
\multicolumn{1}{|l|}{}& \dwhybrid &    3600 & 11.8 & 0.6 & 12 & 0 & 15     &   6879 & 9 & 0.04 & 73 	 \\
\hline
\multicolumn{1}{|l|}{\multirow{5}{*}{\shortstack{\texttt{CloudLarge}\\ ($|\cN|=1000$)}}} & \bsocpcomp  &   3600 & 100.0 & 0.0 & 0 & 0 & 0 &   - & - & - & - 	 \\
\multicolumn{1}{|l|}{}& \dwbc &    3600 & 59.6 & 0.0 & 2 & 0 & 0 &   741 & 100 & 0.04 & 89     \\
 \multicolumn{1}{|l|}{}& \dwpwl &   3600 & 43.1 & 0.2 & 1 & 0 & 6  &  2105 & 100 & 0.01 & 63 		 \\
\multicolumn{1}{|l|}{}& \dwhybrid &    3600 & 34.2 & 0.4 & 1 & 0 & 11 &   4257 & 4 & 0.01 & 8 		 \\
\hline     
\end{tabular}}
\caption{Aggregated statistics of the main computational results} \label{tab.ag}
\end{table}

We first compare the compact BSOCP formulation of SMBP \eqref{master.submod} with the set cover formulation \eqref{master.dw}. So, we evaluate the performance of \bsocpcomp and \dwbc. For all the benchmarks, \dwbc achieves smaller dual gaps than  \bsocpcomp. For small instances, \dwbc also  explores a smaller number of nodes of the search tree. These observations agree with \Cref{prop.strong}  that the continuous relaxation of the set covering formulation is stronger than the continuous relaxation of the compact formulation. The number of nonlinear integer constraints in the compact formulation increases with the number of bins. For medium instances,   \bsocpcomp cannot even finish the root node computation of the compact formulation. However, \dwbc can prove a dual gap or improve primal solutions by solving the set cover formulation. Although the two formulations are insufficient to tackle medium or large instances, the compact formulation is better overall than the set cover formulation.  Then, we will solely examine algorithms that tackle the set cover formulation in the following.

We next evaluate our core innovation to solve the pricing subproblems: the PWL relaxation and its associated combined formulation \eqref{pricing.rel.comb}. So, we compare  \dwbc with  \dwpwl. \dwbc just calls \texttt{CPLEX} to solve the BSCOP formulation \eqref{pricing.ref.bsocp} of pricing subproblems, while \dwpwl uses a tailored branch-and-cut algorithm to solve the combined formulation \eqref{pricing.rel.comb}. Looking at  the problem statistics for all the benchmarks, we find that \dwpwl significantly reduces the master problem' dual gap than \dwbc. Especially for small instances, \dwpwl achieves nearly five times improvement to \dwbc. More details can be found in pricing statistics. \dwpwl can solve pricing subproblems to optimality (pricing gap on average is $0.01\%$) in a short time and thus produce much more columns than \dwbc. Especially for large instances, we find that combined formulation \eqref{pricing.rel.comb} is still solvable. The overall quality of  the combined formulation for submodular knapsack outperforms that of the BSOCP formulation \eqref{pricing.ref.bsocp}.

We examine the hybrid pricing strategy, which replaces the computationally expensive exact pricing with computationally cheap heuristic pricing when the exact pricing is not in need. So, we compare \dwpwl with \dwhybrid. Looking at the problem statistics, the hybrid pricing strategy achieves smaller dual gaps, especially for large instances; it also saves computational time for small instances.  Looking at the pricing statistics, the hybrid pricing strategy can generate twice the number of columns than the exact pricing, for all the instances. As a byproduct, with more columns, \texttt{SCIP} can find more improved primal solutions. We find that the hybrid pricing strategy gives rise to consistent improvement.

We look at the column selection heuristic. So, we compare  \dwhybrid with \dwhybrids. The column selection heuristic can find more improved solutions.

Finally, we summarize our computational results. The set cover formulation is better than the compact formulation regarding scalability, although both formulations are unsolvable for medium and large instances. Our techniques can improve the column generation procedure  for the set cover formulation. Regarding pricing subproblems, a dense BSOCP constraint might be reformulated as a submodular knapsack constraint, so the good performance of PWL relaxations  suggests that PWL relaxations can provide strong MILP relaxations for dense BSOCP constraints. This finding can also help solve other BSOCP problems. The hybrid pricing strategy uses a hint from the Farley bound, so it reduces computational time and is applicable for other column generation problems. As for benchmarks, \texttt{CloudSmall} is a suitable testbed for comparing solvers, \texttt{CloudMedium} is suitable for testing the pricing algorithms, and \texttt{CloudLarge} is still too big to handle.

\section{Conclusion}
\label{sec.conc}

 We develop a PWL-B$\&$C algorithm for solving pricing submodular knapsack problems.
The PWL-B$\&$C algorithm is more efficient than the conventional LP-B$\&$C algorithm implemented in \texttt{CPLEX} for the pricing submodular knapsack problems. The PWL-B$\&$C algorithm can also be extended to solve the multiple submodular knapsack problems. In general MINLP problems, if a nonlinear constraint can be reformulated into a linear part and a univariate concave part, then the univariate concave part can be convexified by PWL relaxation.

Our hybrid pricing strategy applies to the column generation procedure, where the master problems are in set cover formulations, as long as there are fast pricing heuristics. This pricing strategy is helpful for large instances. As a future study, we can apply this strategy to solve the DW decomposition of the capacitated vehicle routing problem, for which the pricing problem is complex.
 
The primary efforts of this paper are solving pricing submodular knapsack subproblems with conflicts  via PWL relaxations and speeding up column generation via a hybrid pricing strategy. There is still much room for improvement in future studies.  Since the submodular knapsack with conflicts is solved multiple times with different parameters, the information of previous column generation iterations can be leveraged statistically to reduce the search space of pricing subproblems.

On the other hand,  commonly known techniques for branch-and-price algorithms are generally helpful. Combining our techniques with other advanced elements from general-purpose framework \cite{pessoa2021solving} could be also useful. For example, we can use stabilization techniques to speed up the convergence of the column generation or use cutting planes to tighten the relaxation of the master problem.

 \section*{ACKNOWLEDGMENT}%
The authors would like to thank Leo Liberti and Sandra Ulrich Ngueveu for discussing this paper with the authors.

\appendix

\section{Proof of \Cref{prop.strong}}
\label{proof.prop.strong}
 \begin{proof}
 Let $$F_j := \{(v_{1j},\cdots,v_{nj}, y_j) \in \{0,1\}^{n+1} : \sum_{i \in \cN} a_i v_{ij} + \sigma\sqrt{\sum_{i \in \cN} b_i v_{ij} } \le c  {y}_{{j}}\}$$ be the feasible set of the $j-$th constraint in the BSOCP formulation of \eqref{master.submod}.  Therefore, the feasible set of the BSCOP formulation is $F = \prod_{j \in \cM} F_j$.
 
Let $\overline{F}_j$  be the continuous relaxation of $F_j$, and
$$\overline{F}_j = \{(v_{1j},\cdots,v_{nj}, y_j) \in [0,1]^{n+1}: \sum_{i \in \cN} a_i v_{ij} + \sigma\sqrt{\sum_{i \in \cN} b_i v_{ij} } \le c  {y}_{{j}}\}.$$ Therefore, the feasible set of the continuous relaxation of the BSCOP formulation is $\overline{F}=\prod_{j \in \cM} \overline{F}_j$.

 On the other hand, the points of $F_j$ are zero vectors and $(p,1)$ ($p \in \cP$). Therefore, its convex hull is
 \begin{multline*}
	 \mathrm{conv}(F_j) = \\
	 \{(v_{1j},\cdots,v_{nj},y_j)  \in [0,1]^{n+1}: \exists \lambda_p \in [0,1]^{\cP} \land \sum_{p \in \cP} \lambda_p = y_j \land v =  \sum_{p \in \cP} d_p \lambda_p \}.
 \end{multline*}
 We note that $\overline{F}_j$ is also a convex relaxation of $F_j$, hence $F_j  \subset \mathrm{conv}(F_j) \subset \overline{F}_j$.

 The optimum of the continuous relaxation of the BSOCP formulation is $\displaystyle \min_{(v,y) \in \overline{F}, v \textup{ satisfies \eqref{master.submod.cover}}} \sum_{j \in \cM} y_j$.
 
  An optimal solution of the LP relaxation of the set cover formulation satisfies $\sum_{p \in \cP} d_{ip} \lambda_p = 1$ ($ i \in \cN$), and the optimal value is exactly the same as $\displaystyle \min_{(v,y) \in \prod_{j \in \cM} \mathrm{conv}(F_j), v \textup{ satisfies \eqref{master.submod.cover}}} \sum_{j \in \cM} y_j$. Since $\prod_{j \in \cM} \mathrm{conv}(F_j) \subset \overline{F}$, the result follows.
   \end{proof}

\section{Proof of \Cref{lem.approx}}
\label{append.lem.approx}
\begin{proof}
Since $\bar{q}_{\cB}$ and $q$ have the same value at $ w \in \{w_1,\hdots,w_{h}\}$, it follows that the $\ell_\infty$-norm is the maximum value of $\ell_\infty$-norms over individual sub intervals:
\begin{equation*}
\ell_\infty(\bar{q}_{\cB}, q) = \max_{w \in [\underline{w},\overline{w}] }|\bar{q}_{\cB}(w) - q(w)| =  \max_{2\le k \le h} \max_{w \in [w_{k-1}, w_k]} |\bar{q}_{\cB}(w) - q(w)|.
\end{equation*}

Let $w \in [w_{k-1}, w_k]$, then
\begin{equation*}
\begin{split}
 & |\bar{q}_{\cB}(w) - q(w)| \\
  =&\frac{q(w_k)-q(w_{k-1})}{w_{k}-w_{k-1}}(w-w_{k-1}) + q(w_{k-1}) - (c-w)^2\\
  =& (w - w_{k-1})(w_k - w).
\end{split}
\end{equation*}
We have
\begin{equation*}
\begin{split}
 &  \max_{w \in [w_{k-1}, w_k]} |\bar{q}_{\cB}(w) - q(w)| \\
  =&  \max_{w \in [w_{k-1}, w_k]}  (w - w_{k-1})(w_k - w)\\
  =& \frac{(w_k -w_{k-1})^2}{4}.
\end{split}
\end{equation*}
The maximum value is at $w =  \frac{w_{k-1}+w_k}{2} $.

It follows that \eqref{bsp} is equivalent to:
\begin{equation*}
	 \min_{\underline{w} = w_1 \le \hdots \le w_h = \overline{w}} \max_{2\le k \le h}\frac{(w_k -w_{k-1})^2}{4}.
\end{equation*}
Therefore, the optimal solution is an equidistant  partition of $[\underline{w}, \overline{w}]$, and the results follow.
 \end{proof}

\section{Non-equidistant breakpoints}
\label{sec.effect}

According to \Cref{lem.approx} in \Cref{sec.pwlrel}, the optimal breakpoints under the $\ell_\infty$ error  form an equidistant partition of $[\underline{w}, \overline{w}]$. In this section, we investigate whether adaptive non-equidistant breakpoints can improve the \dwpwl. There are many possibilities for non-equidistant breakpoints, and we propose a regression approach using the previous pricing information.

We recall that the \dwpwl solver adds a lazy cut at  each infeasible solution $\hat{x}$. Let $\hat{w} := \sum_{i \in \cN'} a'_i \hat{x}_i$ be the corresponding value of variable $w$, and we call it an infeasible $w$-value. For an objective coefficient vector $c'$, let $[w_{\ell}(c'),w_{u}(c')]$ be the range of the set of infeasible $w$-values, which are recorded  during the PWL-B$\&$C algorithm for every pricing problem. Our intuition is that for a new pricing problem with an  objective coefficient vector $c$,  one may reduce the search space by concentrating breakpoints  to the range $[w_{\ell}(c),w_{u}(c)]$, because this  refines the PWL relaxation in that region. Usually $[w_{\ell}(c),w_{u}(c)]$ is unknown, so one can only use a predication range $[w'_{\ell}(c), w'_{u}(c)]$. Given the fixed number of breakpoints of $\cB$, with this limited resource, we use the $k$nn regression approach to learn $[w'_{\ell}(c), w'_{u}(c)]$ and  concentrate a subset of of $\cB$ to $[w'_{\ell}(c), w'_{u}(c)]$.

The  $k$nn regression is as follows. Let $T$ be the number of pricing iterations, and we use a list $\{[w_{\ell}(c^t), w_u(c^t)]\}_{1 \le t \le T}$ to record the set  of intervals, where $w_{\ell}(c^t), w_u(c^t)$ are the lower and upper bounds of the set of infeasible $w$-values in the $t$-th pricing problem. For the new objective coefficient vector $c$, we sort the list in an increasing order w.r.t. the $\ell_2$-norm distances between  $\{c^t\}_{1 \le t \le T}$ to $c$.  The predicted range $[w'_{\ell}(c), w'_{u}(c)]$ is as follows:
\begin{equation*}
    w'_{\ell}(c)= \sum_{1 \le t \le k} w_{\ell}(c^t)  /k,  w'_{u}(c)= \sum_{1 \le t \le k} w_{u}(c^t)  /k.
\end{equation*}
  Recall that there are in total \(h\) breakpoints  in the range \([\underline{w}, \overline{w}]\). Let $r:=(w'_{u}(c) -  w'_{\ell}(c)) / (\overline{w} - \underline{w})$ be the range ratio.
Then, given a concentration scale $ s > 1$,
we put $h* r * s$ number of breakpoints equidistantly in $[w'_{\ell}(c), w'_{u}(c)]$, and $h* (1 - r * s)$ number of breakpoints equidistantly in the remaining breakpoint region. This concentration results in a PWL relaxation with a better approximation in  $[w'_{\ell}(c), w'_{u}(c)]$. Therefore, we hope that the adaptive PWL relaxation could use the previous pricing information.

To understand the performance of the  $k$nn regression approach, we have several configurations with combinations of $k \in \{1,3,5\}$ and $s \in \{ 1.5, 2, 2.5\}$. We note that with $k=1,s = 1$, the configuration is exactly the \dwpwl solver. To test these configurations, we generate two new benchmarks with $|\cN| \in \{500, 900\}$, each containing 36 instances. The aggregated computational results are presented in \Cref{tab.adapt}, which displays SGMs of the relative dual gap of master problems, and the number of generated columns.

 The non-equidistant breakpoints  do not lead to an improvement of the algorithm. In most cases, $k$nn regression approach is even worse than  the equidistant breakpoint approach. Therefore, finding good breakpoints is a complex task.
 
\begin{table}[]
\centering
\resizebox{0.99\columnwidth}{!}{
\renewcommand{\arraystretch}{1.7}
\begin{tabular}{|l*{4}{r}|*{3}{r}|*{3}{r}|}
 \hline
  \multirow{2}{*}{}  &     \multicolumn{4}{|c}{$k=1$} & \multicolumn{3}{|c}{$k=3$} &  \multicolumn{3}{|c|}{$k=5$} \\ [-0.1ex]
   	 &  \multicolumn{1}{|c}{$s=1$}   &  $s=1.5$ & $s=1.5$ &  $s=2.5$  &  $s=1.5$  &  $s=2$ &  $s=2.5$ &  $s=1.5$  &  $s=2$ &  $s=2.5$     \\ [-0.1ex]
		 \cline{1-11}
   	 \multicolumn{1}{|c|}{$\delta_d \%$} &  35.29&35.58&35.66&35.36&35.28&36.28&35.29&35.39&34.95&35.29 \\
   	 \multicolumn{1}{|c|}{$\#C$} & 2134.12&2120.03&2097.46&2123.38&2118.57&2091.97&2148.48&2111.01&2155.95&2134.12 \\
   	 \hline
\end{tabular}}
\caption{Master and pricing problem statistics of different configurations} \label{tab.adapt}
\end{table}

\section{Detailed results}
\label{append.detail}

 Master problem statistics are summarized in \Cref{tab.master.cloudsmall}, \Cref{tab.master.cloudmed}, and \Cref{tab.master.CloudLarge}. Pricing problem statistics are summarized in \Cref{tab.price.cloudsmall}, \Cref{tab.price.cloudmed}, and \Cref{tab.price.CloudLarge}.

 As for the case notation, ``G" denotes the instances of the Gaussian distribution case,  ``H" denotes the instances of the Hoeffding inequality case, and ``D" denotes the instances of the distributionally robust case.

 For each benchmark, we report the SGM statistics of performance metrics, the number of solved instances (\#S), and the number of instances with improved primal bounds  (\#I).

 We also divide the instances in each benchmark into  small subsets and report SGM statistics of these subsets. In these subsets, instances have the same risk level and case.
 
\begin{landscape}
\begin{table}[]
\resizebox{\linewidth}{!}{
\renewcommand{\arraystretch}{1.7}
\begin{tabular}{|l|r|*{6}{r}|*{6}{r}|*{6}{r}|*{6}{r}| *{6}{r}|}
 \hline
		 \multirow{2}{*}{\Shortstack{C a s e}}  &     \multirow{2}{*}{$\alpha$}  &    
		 \multicolumn{6}{|l}{\bsocpcomp} & \multicolumn{6}{|l}{\dwbc} & \multicolumn{6}{|l}{\dwpwl} & \multicolumn{6}{|l}{\dwhybrid}& \multicolumn{6}{|l|}{\dwhybrids} \\ \cline{3-32}
		 & & $t$ & $\delta_d\%$ & $\delta_p\%$ & {\#N} & {\#S} & {\#I}  & $t$ & $\delta_d\%$ & $\delta_p\%$ & {\#N} & {\#S} & {\#I} & $t$ & $\delta_d\%$ & $\delta_p\%$ & {\#N} & {\#S} & {\#I} & $t$ & $\delta_d\%$ & $\delta_p\%$ & {\#N} & {\#S} & {\#I} & $t$ & $\delta_d\%$ & $\delta_p\%$ & {\#N} & {\#S} & {\#I}
		 \\  \hline
\multirow{6}{*}{G}&0.6 & 38 & 1.6 & 0.0 & 70 & 4 & 0 & 211 & 1.7 & 0.0 & 7 & 4 & 0 & 43 & 0.0 & 1.2 & 2 & 6 & 1 & 47 & 0.6 & 0.0 & 7 & 5 & 0 & 42 & 0.6 & 0.0 & 6 & 5 & 0 \\&0.7 & 1199 & 10.6 & 0.0 & 42897 & 1 & 0 & 1049 & 6.0 & 0.0 & 10 & 2 & 0 & 145 & 1.6 & 0.0 & 9 & 4 & 0 & 80 & 1.5 & 0.0 & 20 & 4 & 0 & 84 & 1.5 & 0.0 & 19 & 4 & 0 \\&0.8 & 3600 & 21.6 & 0.0 & 215412 & 0 & 0 & 3600 & 19.3 & 0.0 & 21 & 0 & 0 & 3161 & 10.0 & 0.0 & 287 & 1 & 0 & 2792 & 9.7 & 0.0 & 1177 & 1 & 0 & 2917 & 9.8 & 0.0 & 918 & 1 & 0 \\&0.9 & 3600 & 24.6 & 0.0 & 236150 & 0 & 0 & 3600 & 23.3 & 0.0 & 25 & 0 & 0 & 624 & 1.7 & 3.7 & 79 & 4 & 2 & 221 & 1.5 & 3.7 & 112 & 4 & 2 & 271 & 1.5 & 3.7 & 126 & 4 & 2 \\&0.95 & 3600 & 30.3 & 0.0 & 125945 & 0 & 0 & 3600 & 23.3 & 0.0 & 16 & 0 & 0 & 1160 & 3.6 & 3.7 & 140 & 3 & 2 & 1190 & 5.1 & 6.5 & 743 & 2 & 3 & 836 & 1.4 & 34.1 & 465 & 4 & 5 \\&0.99 & 3600 & 35.6 & 0.0 & 85454 & 0 & 0 & 3600 & 27.4 & 0.0 & 25 & 0 & 0 & 1952 & 5.5 & 6.3 & 266 & 2 & 3 & 1624 & 10.1 & 0.8 & 1110 & 1 & 1 & 828 & 5.4 & 6.3 & 442 & 2 & 3 \\\hline\multirow{6}{*}{H}&0.6 & 535 & 3.5 & 0.0 & 6828 & 3 & 0 & 2300 & 6.4 & 1.2 & 66 & 2 & 1 & 777 & 1.7 & 1.2 & 35 & 4 & 1 & 166 & 0.6 & 3.7 & 41 & 5 & 2 & 253 & 1.7 & 1.2 & 43 & 4 & 1 \\&0.7 & 532 & 6.0 & 0.0 & 3565 & 2 & 0 & 572 & 3.1 & 0.0 & 6 & 3 & 0 & 71 & 0.0 & 0.0 & 1 & 6 & 0 & 18 & 0.0 & 0.0 & 1 & 6 & 0 & 19 & 0.0 & 0.0 & 1 & 6 & 0 \\&0.8 & 178 & 3.3 & 0.0 & 2159 & 3 & 0 & 1728 & 7.6 & 0.0 & 37 & 2 & 0 & 572 & 0.7 & 9.0 & 83 & 5 & 3 & 528 & 3.8 & 1.2 & 277 & 3 & 1 & 587 & 3.8 & 1.2 & 247 & 3 & 1 \\&0.9 & 352 & 3.4 & 0.0 & 1807 & 3 & 0 & 3600 & 20.4 & 0.0 & 76 & 0 & 0 & 1506 & 6.1 & 0.0 & 160 & 2 & 0 & 413 & 1.6 & 3.7 & 167 & 4 & 2 & 455 & 1.6 & 3.7 & 162 & 4 & 2 \\&0.95 & 734 & 6.5 & 0.0 & 23189 & 2 & 0 & 990 & 6.5 & 0.0 & 17 & 2 & 0 & 329 & 1.5 & 0.0 & 18 & 4 & 0 & 119 & 0.6 & 1.2 & 17 & 5 & 1 & 203 & 1.4 & 1.2 & 27 & 4 & 1 \\&0.99 & 3600 & 20.9 & 0.0 & 345213 & 0 & 0 & 3600 & 19.7 & 0.0 & 16 & 0 & 0 & 253 & 1.8 & 1.2 & 31 & 4 & 1 & 215 & 1.5 & 1.2 & 91 & 4 & 1 & 222 & 2.9 & 0.9 & 81 & 3 & 1 \\\hline\multirow{6}{*}{D}&0.6 & 3600 & 24.0 & 0.0 & 187687 & 0 & 0 & 3600 & 20.2 & 0.0 & 18 & 0 & 0 & 419 & 1.7 & 1.2 & 34 & 4 & 1 & 145 & 0.6 & 7.8 & 55 & 5 & 3 & 135 & 1.5 & 3.1 & 46 & 4 & 2 \\&0.7 & 3600 & 29.0 & 0.0 & 153161 & 0 & 0 & 3600 & 23.1 & 2.3 & 21 & 0 & 2 & 1234 & 1.6 & 9.0 & 141 & 4 & 3 & 458 & 1.6 & 9.0 & 276 & 4 & 3 & 476 & 1.6 & 9.0 & 246 & 4 & 3 \\&0.8 & 3600 & 34.7 & 0.0 & 95858 & 0 & 0 & 3600 & 22.3 & 0.0 & 20 & 0 & 0 & 3600 & 14.4 & 0.9 & 463 & 0 & 1 & 2153 & 2.8 & 9.0 & 1542 & 3 & 3 & 2476 & 4.6 & 7.8 & 1305 & 2 & 3 \\&0.9 & 3600 & 44.6 & 0.0 & 41259 & 0 & 0 & 3600 & 25.8 & 2.4 & 27 & 0 & 2 & 3146 & 8.8 & 6.6 & 384 & 1 & 3 & 2389 & 9.0 & 3.0 & 1260 & 1 & 2 & 2459 & 7.4 & 25.0 & 1111 & 1 & 5 \\&0.95 & 3600 & 62.1 & 0.0 & 29645 & 0 & 0 & 3600 & 15.4 & 46.9 & 20 & 0 & 6 & 3321 & 6.3 & 31.8 & 414 & 1 & 5 & 2332 & 3.5 & 66.3 & 664 & 2 & 6 & 2028 & 3.8 & 59.1 & 741 & 2 & 6 \\&0.99 & 3600 & 77.5 & 0.0 & 10352 & 0 & 0 & 1218 & 0.6 & 96.4 & 33 & 5 & 6 & 117 & 0.0 & 100.0 & 90 & 6 & 6 & 30 & 0.0 & 100.0 & 32 & 6 & 6 & 27 & 0.0 & 100.0 & 32 & 6 & 6 \\\hline\multicolumn{2}{|c|}{All}&1452 & 15.8 & 0.0 & 26601 & 18 & 0 & 2129 & 11.4 & 0.9 & 21 & 20 & 17 & 633 & 2.4 & 2.7 & 66 & 61 & 32 & 330 & 2.0 & 3.4 & 127 & 65 & 36 & 335 & 2.1 & 4.1 & 114 & 63 & 41 \\ \hline
\end{tabular}}
\caption{Master problem statistics of \texttt{CloudSmall} with 108 instances ($|\cN|=100$)} \label{tab.master.cloudsmall}
\end{table}

\begin{table}[]
\resizebox{\linewidth}{!}{
\renewcommand{\arraystretch}{1.7}
\begin{tabular}{|l|r|*{6}{r}|*{6}{r}|*{6}{r}|*{6}{r}| *{6}{r}|}
 \hline
		 \multirow{2}{*}{\Shortstack{C a s e}}  &     \multirow{2}{*}{$\alpha$}  &    
		 \multicolumn{6}{|l}{\bsocpcomp} & \multicolumn{6}{|l}{\dwbc} & \multicolumn{6}{|l}{\dwpwl} & \multicolumn{6}{|l}{\dwhybrid}& \multicolumn{6}{|l|}{\dwhybrids} \\ \cline{3-32}
		 & & $t$ & $\delta_d\%$ & $\delta_p\%$ & {\#N} & {\#S} & {\#I}  & $t$ & $\delta_d\%$ & $\delta_p\%$ & {\#N} & {\#S} & {\#I} & $t$ & $\delta_d\%$ & $\delta_p\%$ & {\#N} & {\#S} & {\#I} & $t$ & $\delta_d\%$ & $\delta_p\%$ & {\#N} & {\#S} & {\#I} & $t$ & $\delta_d\%$ & $\delta_p\%$ & {\#N} & {\#S} & {\#I}
		 \\  \hline
\multirow{6}{*}{G}&0.6 & 3600 & 100.0 & 0.0 & 0 & 0 & 0 & 3600 & 20.6 & 0.0 & 2 & 0 & 0 & 3600 & 14.1 & 0.0 & 2 & 0 & 0 & 3600 & 8.5 & 0.0 & 32 & 0 & 0 & 3600 & 8.6 & 0.0 & 32 & 0 & 0 \\&0.7 & 3600 & 100.0 & 0.0 & 0 & 0 & 0 & 3600 & 30.1 & 0.0 & 2 & 0 & 0 & 3600 & 16.3 & 0.0 & 2 & 0 & 0 & 3600 & 10.9 & 0.0 & 122 & 0 & 0 & 3600 & 10.9 & 0.0 & 100 & 0 & 0 \\&0.8 & 3600 & 100.0 & 0.0 & 0 & 0 & 0 & 3600 & 36.4 & 0.0 & 2 & 0 & 0 & 3600 & 15.9 & 0.0 & 2 & 0 & 0 & 3600 & 11.3 & 0.0 & 46 & 0 & 0 & 3600 & 9.7 & 4.4 & 42 & 0 & 3 \\&0.9 & 3600 & 100.0 & 0.0 & 0 & 0 & 0 & 3600 & 42.2 & 0.0 & 2 & 0 & 0 & 3600 & 19.1 & 0.0 & 2 & 0 & 0 & 3600 & 13.0 & 0.0 & 18 & 0 & 0 & 3600 & 11.2 & 11.0 & 24 & 0 & 5 \\&0.95 & 3600 & 100.0 & 0.0 & 0 & 0 & 0 & 3600 & 48.4 & 0.0 & 2 & 0 & 0 & 3600 & 18.8 & 0.0 & 2 & 0 & 0 & 3600 & 14.8 & 0.8 & 9 & 0 & 1 & 3600 & 13.7 & 14.8 & 8 & 0 & 6 \\&0.99 & 3600 & 100.0 & 0.0 & 0 & 0 & 0 & 3600 & 47.6 & 0.0 & 2 & 0 & 0 & 3600 & 22.3 & 0.0 & 2 & 0 & 0 & 3600 & 14.1 & 0.0 & 2 & 0 & 0 & 3600 & 13.3 & 2.9 & 2 & 0 & 3 \\\hline\multirow{6}{*}{H}&0.6 & 3600 & 100.0 & 0.0 & 0 & 0 & 0 & 3600 & 27.7 & 0.0 & 2 & 0 & 0 & 3600 & 16.4 & 0.0 & 2 & 0 & 0 & 3600 & 10.5 & 0.0 & 22 & 0 & 0 & 3600 & 10.5 & 0.0 & 21 & 0 & 0 \\&0.7 & 3600 & 100.0 & 0.0 & 0 & 0 & 0 & 3600 & 29.0 & 0.0 & 2 & 0 & 0 & 3600 & 18.0 & 0.0 & 2 & 0 & 0 & 3600 & 10.5 & 0.0 & 30 & 0 & 0 & 3600 & 10.6 & 0.0 & 32 & 0 & 0 \\&0.8 & 3600 & 100.0 & 0.0 & 0 & 0 & 0 & 3600 & 31.5 & 0.0 & 2 & 0 & 0 & 3600 & 17.2 & 0.0 & 2 & 0 & 0 & 3600 & 11.6 & 0.0 & 9 & 0 & 0 & 3600 & 11.6 & 0.0 & 10 & 0 & 0 \\&0.9 & 3600 & 100.0 & 0.0 & 0 & 0 & 0 & 3600 & 33.3 & 0.0 & 2 & 0 & 0 & 3600 & 17.0 & 0.0 & 2 & 0 & 0 & 3600 & 10.9 & 0.0 & 7 & 0 & 0 & 3600 & 10.8 & 0.0 & 9 & 0 & 0 \\&0.95 & 3600 & 100.0 & 0.0 & 0 & 0 & 0 & 3600 & 35.1 & 0.0 & 2 & 0 & 0 & 3600 & 17.6 & 0.0 & 2 & 0 & 0 & 3600 & 10.7 & 0.0 & 24 & 0 & 0 & 3600 & 10.7 & 0.0 & 23 & 0 & 0 \\&0.99 & 3600 & 100.0 & 0.0 & 0 & 0 & 0 & 3600 & 39.3 & 0.0 & 2 & 0 & 0 & 3600 & 19.4 & 0.0 & 2 & 0 & 0 & 3600 & 12.7 & 0.0 & 27 & 0 & 0 & 3600 & 12.7 & 0.0 & 27 & 0 & 0 \\\hline\multirow{6}{*}{D}&0.6 & 3600 & 100.0 & 0.0 & 0 & 0 & 0 & 3600 & 42.9 & 0.0 & 2 & 0 & 0 & 3600 & 20.2 & 0.9 & 2 & 0 & 1 & 3600 & 12.9 & 0.0 & 10 & 0 & 0 & 3600 & 10.3 & 21.8 & 13 & 0 & 6 \\&0.7 & 3600 & 100.0 & 0.0 & 0 & 0 & 0 & 3600 & 48.7 & 0.0 & 2 & 0 & 0 & 3600 & 21.3 & 0.0 & 2 & 0 & 0 & 3600 & 15.3 & 0.6 & 12 & 0 & 1 & 3600 & 13.6 & 15.0 & 19 & 0 & 6 \\&0.8 & 3600 & 100.0 & 0.0 & 0 & 0 & 0 & 3600 & 48.4 & 0.0 & 2 & 0 & 0 & 3600 & 22.7 & 0.0 & 2 & 0 & 0 & 3600 & 15.1 & 0.9 & 2 & 0 & 1 & 3600 & 14.8 & 7.4 & 2 & 0 & 5 \\&0.9 & 3600 & 100.0 & 0.0 & 0 & 0 & 0 & 3600 & 51.2 & 0.0 & 2 & 0 & 0 & 3600 & 21.0 & 0.0 & 2 & 0 & 0 & 3600 & 15.4 & 0.0 & 2 & 0 & 0 & 3600 & 14.8 & 2.3 & 2 & 0 & 3 \\&0.95 & 3600 & 100.0 & 0.0 & 0 & 0 & 0 & 3600 & 55.9 & 0.0 & 2 & 0 & 0 & 3600 & 22.3 & 2.6 & 2 & 0 & 3 & 3600 & 17.4 & 13.7 & 2 & 0 & 6 & 3600 & 16.5 & 17.6 & 2 & 0 & 6 \\&0.99 & 3600 & 100.0 & 0.0 & 0 & 0 & 0 & 3600 & 57.4 & 3.3 & 2 & 0 & 4 & 3600 & 4.0 & 78.4 & 2 & 0 & 6 & 3600 & 4.2 & 76.4 & 3 & 0 & 6 & 3600 & 3.9 & 77.4 & 2 & 0 & 6 \\\hline\multicolumn{2}{|c|}{All}&3600 & 100.0 & 0.0 & 0 & 0 & 0 & 3600 & 39.0 & 0.1 & 2 & 0 & 4 & 3600 & 17.2 & 0.4 & 1 & 0 & 10 & 3600 & 11.8 & 0.6 & 12 & 0 & 15 & 3600 & 11.2 & 3.0 & 12 & 0 & 49 \\\hline
\end{tabular}}
\caption{Master problem statistics of \texttt{CloudMedium} with 108 instances ($|\cN|=400$)} \label{tab.master.cloudmed}
\end{table}

\begin{table}[]
\resizebox{\linewidth}{!}{
\renewcommand{\arraystretch}{1.7}
\begin{tabular}{|l|r|*{6}{r}|*{6}{r}|*{6}{r}|*{6}{r}| *{6}{r}|}
 \hline
		 \multirow{2}{*}{\Shortstack{C a s e}}  &     \multirow{2}{*}{$\alpha$}  &    
		 \multicolumn{6}{|l}{\bsocpcomp} & \multicolumn{6}{|l}{\dwbc} & \multicolumn{6}{|l}{\dwpwl} & \multicolumn{6}{|l}{\dwhybrid}& \multicolumn{6}{|l|}{\dwhybrids} \\ \cline{3-32}
		 & & $t$ & $\delta_d\%$ & $\delta_p\%$ & {\#N} & {\#S} & {\#I}  & $t$ & $\delta_d\%$ & $\delta_p\%$ & {\#N} & {\#S} & {\#I} & $t$ & $\delta_d\%$ & $\delta_p\%$ & {\#N} & {\#S} & {\#I} & $t$ & $\delta_d\%$ & $\delta_p\%$ & {\#N} & {\#S} & {\#I} & $t$ & $\delta_d\%$ & $\delta_p\%$ & {\#N} & {\#S} & {\#I}
		 \\   \hline
\multirow{6}{*}{G}&0.6 & 3600 & 100.0 & 0.0 & 0 & 0 & 0 & 3600 & 38.9 & 0.0 & 2 & 0 & 0 & 3600 & 38.7 & 0.0 & 2 & 0 & 0 & 3600 & 32.6 & 0.0 & 2 & 0 & 0 & 3600 & 32.6 & 0.0 & 2 & 0 & 0 \\&0.7 & 3600 & 100.0 & 0.0 & 0 & 0 & 0 & 3600 & 41.7 & 0.0 & 2 & 0 & 0 & 3600 & 38.7 & 0.0 & 2 & 0 & 0 & 3600 & 34.9 & 0.0 & 2 & 0 & 0 & 3600 & 35.2 & 0.0 & 2 & 0 & 0 \\&0.8 & 3600 & 100.0 & 0.0 & 0 & 0 & 0 & 3600 & 51.2 & 0.0 & 2 & 0 & 0 & 3600 & 39.1 & 0.0 & 2 & 0 & 0 & 3600 & 35.0 & 0.0 & 2 & 0 & 0 & 3600 & 35.2 & 0.0 & 2 & 0 & 0 \\&0.9 & 3600 & 100.0 & 0.0 & 0 & 0 & 0 & 3600 & 64.1 & 0.0 & 2 & 0 & 0 & 3600 & 47.3 & 0.0 & 2 & 0 & 0 & 3600 & 41.3 & 0.2 & 1 & 0 & 1 & 3600 & 41.9 & 0.9 & 1 & 0 & 3 \\&0.95 & 3600 & 100.0 & 0.0 & 0 & 0 & 0 & 3600 & 71.3 & 0.0 & 2 & 0 & 0 & 3600 & 47.7 & 0.0 & 1 & 0 & 0 & 3600 & 41.5 & 0.2 & 1 & 0 & 1 & 3600 & 41.6 & 1.2 & 1 & 0 & 4 \\&0.99 & 3600 & 100.0 & 0.0 & 0 & 0 & 0 & 3600 & 79.8 & 0.0 & 2 & 0 & 0 & 3600 & 45.0 & 0.0 & 1 & 0 & 0 & 3600 & 38.4 & 0.2 & 1 & 0 & 1 & 3600 & 38.9 & 0.4 & 1 & 0 & 2 \\\hline\multirow{6}{*}{H}&0.6 & 3600 & 100.0 & 0.0 & 0 & 0 & 0 & 3600 & 44.5 & 0.0 & 2 & 0 & 0 & 3600 & 41.6 & 0.0 & 2 & 0 & 0 & 3600 & 36.4 & 0.0 & 1 & 0 & 0 & 3600 & 36.8 & 0.0 & 2 & 0 & 0 \\&0.7 & 3600 & 100.0 & 0.0 & 0 & 0 & 0 & 3600 & 45.7 & 0.0 & 2 & 0 & 0 & 3600 & 41.0 & 0.0 & 2 & 0 & 0 & 3600 & 35.5 & 0.0 & 2 & 0 & 0 & 3600 & 35.7 & 0.0 & 2 & 0 & 0 \\&0.8 & 3600 & 100.0 & 0.0 & 0 & 0 & 0 & 3600 & 47.5 & 0.0 & 2 & 0 & 0 & 3600 & 43.1 & 0.0 & 2 & 0 & 0 & 3600 & 37.5 & 0.0 & 2 & 0 & 0 & 3600 & 36.1 & 0.0 & 1 & 0 & 0 \\&0.9 & 3600 & 100.0 & 0.0 & 0 & 0 & 0 & 3600 & 50.3 & 0.0 & 2 & 0 & 0 & 3600 & 41.9 & 0.0 & 2 & 0 & 0 & 3600 & 36.6 & 0.0 & 1 & 0 & 0 & 3600 & 36.7 & 0.0 & 2 & 0 & 0 \\&0.95 & 3600 & 100.0 & 0.0 & 0 & 0 & 0 & 3600 & 52.1 & 0.0 & 2 & 0 & 0 & 3600 & 43.3 & 0.0 & 2 & 0 & 0 & 3600 & 37.3 & 0.0 & 1 & 0 & 0 & 3600 & 37.0 & 0.0 & 1 & 0 & 0 \\&0.99 & 3600 & 100.0 & 0.0 & 0 & 0 & 0 & 3600 & 58.5 & 0.0 & 2 & 0 & 0 & 3600 & 46.1 & 0.0 & 2 & 0 & 0 & 3600 & 41.3 & 0.0 & 1 & 0 & 0 & 3600 & 41.3 & 0.0 & 1 & 0 & 0 \\\hline\multirow{6}{*}{D}&0.6 & 3600 & 100.0 & 0.0 & 0 & 0 & 0 & 3600 & 62.2 & 0.0 & 2 & 0 & 0 & 3600 & 48.0 & 0.0 & 1 & 0 & 0 & 3600 & 41.4 & 0.0 & 1 & 0 & 0 & 3600 & 41.7 & 0.0 & 1 & 0 & 0 \\&0.7 & 3600 & 100.0 & 0.0 & 0 & 0 & 0 & 3600 & 69.6 & 0.0 & 2 & 0 & 0 & 3600 & 46.8 & 0.0 & 1 & 0 & 0 & 3600 & 42.4 & 0.0 & 1 & 0 & 0 & 3600 & 42.2 & 0.5 & 1 & 0 & 2 \\&0.8 & 3600 & 100.0 & 0.0 & 0 & 0 & 0 & 3600 & 77.1 & 0.0 & 2 & 0 & 0 & 3600 & 45.3 & 0.0 & 1 & 0 & 0 & 3600 & 39.6 & 0.0 & 1 & 0 & 0 & 3600 & 39.4 & 0.5 & 1 & 0 & 2 \\&0.9 & 3600 & 100.0 & 0.0 & 0 & 0 & 0 & 3600 & 82.2 & 0.0 & 2 & 0 & 0 & 3600 & 44.1 & 0.0 & 2 & 0 & 0 & 3600 & 35.7 & 0.0 & 1 & 0 & 0 & 3600 & 36.7 & 0.0 & 1 & 0 & 0 \\&0.95 & 3600 & 100.0 & 0.0 & 0 & 0 & 0 & 3600 & 84.5 & 0.0 & 2 & 0 & 0 & 3600 & 49.0 & 0.0 & 2 & 0 & 0 & 3600 & 27.7 & 1.6 & 1 & 0 & 2 & 3600 & 29.8 & 7.0 & 2 & 0 & 6 \\&0.99 & 3600 & 100.0 & 0.0 & 0 & 0 & 0 & 3600 & 84.9 & 0.0 & 2 & 0 & 0 & 3600 & 32.9 & 19.9 & 2 & 0 & 6 & 3600 & 7.9 & 61.0 & 2 & 0 & 6 & 3600 & 13.1 & 49.6 & 2 & 0 & 6 \\\hline\multicolumn{2}{|c|}{All}&3600 & 100.0 & 0.0 & 0 & 0 & 0 & 3600 & 59.6 & 0.0 & 2 & 0 & 0 & 3600 & 43.1 & 0.2 & 1 & 0 & 6 & 3600 & 34.2 & 0.4 & 1 & 0 & 11 & 3600 & 35.3 & 0.6 & 1 & 0 & 25 \\\hline
\end{tabular}}
\caption{Master problem statistics of \texttt{CloudLarge} with 108 instances ($|\cN|=1000$)} \label{tab.master.CloudLarge}
\end{table}
\end{landscape}

\begin{table}[]
\centering
\resizebox{\columnwidth}{!}{
\renewcommand{\arraystretch}{1.7}
\begin{tabular}{|l|r|*{4}{r}|*{4}{r}|*{4}{r}|*{4}{r}|}
 \hline
		 \multirow{2}{*}{Case}  &     \multirow{2}{*}{$\alpha$}  &     \multicolumn{4}{|l}{\dwbc} & \multicolumn{4}{|l}{\dwpwl} & \multicolumn{4}{|l}{\dwhybrid}& \multicolumn{4}{|l|}{\dwhybrids} \\ [-0.8ex]
		 \cline{3-18}
		 & & {\#C} & E\% & {$\tau\%$} &$t_p \%$&  {\#C} & E\% & {$\tau\%$} &$t_p \%$&   {\#C} & E\% & {$\tau\%$} &$t_p \%$&   {\#C} & E\% & {$\tau\%$} &$t_p \%$\\
		 \hline
\multirow{6}{*}{G}&0.6 & 519 & 100 & 0.1 & 99 & 226 & 100 & 0.01 & 99 & 440 & 33 & 0.02 & 97 & 389 & 33 & 0.02 & 94 \\ [-1.1ex]&0.7 & 836 & 100 & 2.1 & 99 & 517 & 100 & 0.01 & 99 & 838 & 27 & 0.02 & 97 & 808 & 28 & 0.02 & 89 \\ [-1.1ex]&0.8 & 1784 & 100 & 5.71 & 99 & 7573 & 100 & 0.01 & 99 & 25935 & 12 & 0.01 & 94 & 21026 & 12 & 0.01 & 70 \\ [-1.1ex]&0.9 & 1711 & 100 & 11.67 & 99 & 1874 & 100 & 0.01 & 99 & 2904 & 16 & 0.01 & 96 & 3260 & 15 & 0.01 & 87 \\ [-1.1ex]&0.95 & 1717 & 100 & 12.59 & 99 & 2887 & 100 & 0.01 & 99 & 14654 & 12 & 0.01 & 95 & 8706 & 13 & 0.01 & 81 \\ [-1.1ex]&0.99 & 1699 & 100 & 16.48 & 99 & 5800 & 100 & 0.01 & 99 & 20729 & 12 & 0.01 & 94 & 10076 & 11 & 0.01 & 73 \\ \hline\multirow{6}{*}{H}&0.6 & 2449 & 100 & 0.29 & 99 & 1869 & 100 & 0.02 & 99 & 1511 & 20 & 0.03 & 96 & 2084 & 19 & 0.02 & 86 \\ [-1.1ex]&0.7 & 523 & 100 & 4.14 & 99 & 204 & 100 & 0.02 & 99 & 179 & 35 & 0.01 & 98 & 194 & 34 & 0.01 & 98 \\ [-1.1ex]&0.8 & 1594 & 100 & 0.29 & 99 & 2149 & 100 & 0.01 & 99 & 7162 & 13 & 0.01 & 94 & 6128 & 13 & 0.01 & 76 \\ [-1.1ex]&0.9 & 2663 & 100 & 0.9 & 99 & 3782 & 100 & 0.01 & 99 & 4476 & 16 & 0.01 & 95 & 4117 & 18 & 0.01 & 86 \\ [-1.1ex]&0.95 & 897 & 100 & 1.48 & 99 & 942 & 100 & 0.02 & 99 & 861 & 32 & 0.03 & 97 & 1293 & 33 & 0.02 & 91 \\ [-1.1ex]&0.99 & 1722 & 100 & 8.05 & 99 & 940 & 100 & 0.01 & 99 & 2420 & 21 & 0.01 & 96 & 2274 & 20 & 0.01 & 81 \\ \hline\multirow{6}{*}{D}&0.6 & 1821 & 100 & 6.4 & 99 & 1306 & 100 & 0.01 & 99 & 2099 & 16 & 0.01 & 95 & 1798 & 16 & 0.01 & 86 \\ [-1.1ex]&0.7 & 1728 & 100 & 12.51 & 99 & 3303 & 100 & 0.01 & 99 & 7434 & 10 & 0.01 & 95 & 6730 & 10 & 0.01 & 81 \\ [-1.1ex]&0.8 & 1697 & 100 & 15.03 & 99 & 7952 & 100 & 0.01 & 99 & 23859 & 12 & 0.01 & 95 & 20313 & 13 & 0.01 & 75 \\ [-1.1ex]&0.9 & 1686 & 100 & 17.62 & 99 & 6862 & 100 & 0.01 & 99 & 20794 & 15 & 0.01 & 96 & 18002 & 15 & 0.01 & 83 \\ [-1.1ex]&0.95 & 1703 & 100 & 16.3 & 99 & 5850 & 100 & 0.01 & 99 & 7948 & 37 & 0.01 & 98 & 8224 & 32 & 0.01 & 95 \\ [-1.1ex]&0.99 & 601 & 100 & 10.74 & 99 & 730 & 100 & 0.01 & 99 & 603 & 23 & 0.01 & 97 & 606 & 22 & 0.01 & 96 \\ \hline\multicolumn{2}{|c|}{All}&1373 & 100 & 3.56 & 99 & 1869 & 100 & 0.01 & 99 & 3485 & 18 & 0.01 & 96 & 3204 & 18 & 0.01 & 84 \\\hline
\end{tabular}}
\caption{Pricing problem statistics of \texttt{CloudSmall} with 108 instances ($|\cN|=100$)} \label{tab.price.cloudsmall}
\end{table}

\begin{table}[]
\centering
\resizebox{0.99\columnwidth}{!}{
\renewcommand{\arraystretch}{1.7}
\begin{tabular}{|l|r|*{4}{r}|*{4}{r}|*{4}{r}|*{4}{r}|}
 \hline
		 \multirow{2}{*}{Case}  &     \multirow{2}{*}{$\alpha$}  &     \multicolumn{4}{|l}{\dwbc} & \multicolumn{4}{|l}{\dwpwl} & \multicolumn{4}{|l}{\dwhybrid}& \multicolumn{4}{|l|}{\dwhybrids} \\ [-0.8ex]
		 \cline{3-18}
		 & & {\#C} & E\% & {$\tau\%$} &$t_p \%$&  {\#C} & E\% & {$\tau\%$} &$t_p \%$&   {\#C} & E\% & {$\tau\%$} &$t_p \%$&   {\#C} & E\% & {$\tau\%$} &$t_p \%$\\
		 \hline
\multirow{6}{*}{G}&0.6 & 2060 & 100 & 0.01 & 94 & 3117 & 100 & 0.01 & 88 & 9368 & 4 & 0.06 & 54 & 8054 & 5 & 0.06 & 55 \\ [-1.1ex]&0.7 & 1413 & 100 & 0.02 & 97 & 3622 & 100 & 0.01 & 87 & 12327 & 4 & 0.02 & 54 & 12105 & 4 & 0.01 & 50 \\ [-1.1ex]&0.8 & 943 & 100 & 0.11 & 98 & 3847 & 100 & 0.01 & 88 & 11031 & 5 & 0.02 & 61 & 9597 & 5 & 0.02 & 55 \\ [-1.1ex]&0.9 & 814 & 100 & 0.36 & 99 & 3277 & 100 & 0.01 & 91 & 6595 & 9 & 0.01 & 74 & 7190 & 9 & 0.01 & 73 \\ [-1.1ex]&0.95 & 658 & 100 & 0.77 & 99 & 3200 & 100 & 0.01 & 93 & 6212 & 11 & 0.01 & 81 & 6150 & 11 & 0.01 & 81 \\ [-1.1ex]&0.99 & 589 & 100 & 2.04 & 99 & 3816 & 100 & 0.01 & 92 & 6622 & 10 & 0.01 & 76 & 6786 & 9 & 0.01 & 74 \\ \hline\multirow{6}{*}{H}&0.6 & 1503 & 100 & 0.02 & 97 & 2985 & 100 & 0.01 & 90 & 6343 & 8 & 0.11 & 70 & 6258 & 8 & 0.1 & 70 \\ [-1.1ex]&0.7 & 1281 & 100 & 0.03 & 97 & 3114 & 100 & 0.01 & 90 & 7692 & 7 & 0.07 & 64 & 7005 & 7 & 0.07 & 64 \\ [-1.1ex]&0.8 & 1246 & 100 & 0.03 & 98 & 3346 & 100 & 0.01 & 91 & 5947 & 9 & 0.09 & 75 & 5715 & 10 & 0.09 & 75 \\ [-1.1ex]&0.9 & 1043 & 100 & 0.06 & 98 & 3409 & 100 & 0.01 & 90 & 6472 & 9 & 0.05 & 73 & 6275 & 9 & 0.05 & 73 \\ [-1.1ex]&0.95 & 867 & 100 & 0.13 & 98 & 3145 & 100 & 0.01 & 91 & 6513 & 9 & 0.02 & 74 & 6570 & 9 & 0.02 & 71 \\ [-1.1ex]&0.99 & 840 & 100 & 0.22 & 99 & 3199 & 100 & 0.01 & 91 & 7448 & 9 & 0.02 & 73 & 7067 & 9 & 0.02 & 72 \\ \hline\multirow{6}{*}{D}&0.6 & 794 & 100 & 0.31 & 99 & 3189 & 100 & 0.01 & 91 & 6691 & 9 & 0.01 & 67 & 6962 & 9 & 0.01 & 70 \\ [-1.1ex]&0.7 & 681 & 100 & 0.57 & 99 & 3191 & 100 & 0.01 & 93 & 6069 & 11 & 0.01 & 80 & 6991 & 10 & 0.01 & 78 \\ [-1.1ex]&0.8 & 609 & 100 & 1.55 & 99 & 3674 & 100 & 0.01 & 93 & 5792 & 11 & 0.01 & 81 & 5848 & 11 & 0.01 & 80 \\ [-1.1ex]&0.9 & 552 & 100 & 2.35 & 99 & 3331 & 100 & 0.01 & 94 & 6703 & 9 & 0.02 & 78 & 6844 & 8 & 0.02 & 78 \\ [-1.1ex]&0.95 & 519 & 100 & 2.7 & 99 & 3129 & 100 & 0.01 & 96 & 5757 & 15 & 0.01 & 90 & 5817 & 14 & 0.01 & 87 \\ [-1.1ex]&0.99 & 448 & 100 & 16.35 & 99 & 4386 & 100 & 0.01 & 97 & 4162 & 18 & 0.12 & 97 & 4126 & 18 & 0.09 & 95 \\ \hline\multicolumn{2}{|c|}{All}&861 & 100 & 0.39 & 98 & 3372 & 100 & 0.01 & 91 & 6879 & 9 & 0.04 & 73 & 6797 & 9 & 0.03 & 71 \\\hline
\end{tabular}}
\caption{Pricing problem statistics of \texttt{CloudMedium} with 108 instances ($|\cN|=400$)} \label{tab.price.cloudmed}
\end{table}

\begin{table}[]
\centering
\resizebox{0.99\columnwidth}{!}{
\renewcommand{\arraystretch}{1.7}
\begin{tabular}{|l|r|*{4}{r}|*{4}{r}|*{4}{r}|*{4}{r}|}
 \hline
		 \multirow{2}{*}{Case}  &     \multirow{2}{*}{$\alpha$}  &     \multicolumn{4}{|l}{\dwbc} & \multicolumn{4}{|l}{\dwpwl} & \multicolumn{4}{|l}{\dwhybrid}& \multicolumn{4}{|l|}{\dwhybrids} \\ [-0.8ex]
		 \cline{3-18}
		 & & {\#C} & E\% & {$\tau\%$} &$t_p \%$&  {\#C} & E\% & {$\tau\%$} &$t_p \%$&   {\#C} & E\% & {$\tau\%$} &$t_p \%$&   {\#C} & E\% & {$\tau\%$} &$t_p \%$\\
		 \hline
\multirow{6}{*}{G}&0.6 & 1786 & 100 & 0.01 & 55 & 1923 & 100 & 0.01 & 47 & 2748 & 6 & 0.01 & 5 & 2766 & 6 & 0.01 & 5 \\ [-1.1ex]&0.7 & 1315 & 100 & 0.01 & 79 & 1903 & 100 & 0.01 & 53 & 2953 & 6 & 0.01 & 5 & 2932 & 6 & 0.01 & 5 \\ [-1.1ex]&0.8 & 1038 & 100 & 0.01 & 90 & 1977 & 100 & 0.01 & 56 & 3182 & 5 & 0.01 & 6 & 3155 & 5 & 0.01 & 6 \\ [-1.1ex]&0.9 & 760 & 100 & 0.02 & 96 & 2130 & 100 & 0.01 & 60 & 4433 & 3 & 0.01 & 5 & 4360 & 3 & 0.01 & 5 \\ [-1.1ex]&0.95 & 629 & 100 & 0.03 & 98 & 2223 & 100 & 0.01 & 61 & 4724 & 3 & 0.01 & 6 & 4597 & 3 & 0.01 & 5 \\ [-1.1ex]&0.99 & 409 & 100 & 0.08 & 99 & 2340 & 100 & 0.01 & 64 & 5005 & 3 & 0.01 & 6 & 4795 & 3 & 0.01 & 6 \\ \hline\multirow{6}{*}{H}&0.6 & 1422 & 100 & 0.01 & 76 & 1624 & 100 & 0.01 & 67 & 2968 & 5 & 0.01 & 10 & 2960 & 5 & 0.01 & 10 \\ [-1.1ex]&0.7 & 1357 & 100 & 0.01 & 80 & 1689 & 100 & 0.01 & 66 & 3058 & 6 & 0.01 & 9 & 3017 & 6 & 0.01 & 9 \\ [-1.1ex]&0.8 & 1242 & 100 & 0.01 & 84 & 1723 & 100 & 0.01 & 67 & 3214 & 5 & 0.01 & 8 & 3174 & 5 & 0.01 & 8 \\ [-1.1ex]&0.9 & 1091 & 100 & 0.01 & 88 & 1694 & 100 & 0.01 & 66 & 3134 & 5 & 0.01 & 9 & 3070 & 5 & 0.01 & 9 \\ [-1.1ex]&0.95 & 1017 & 100 & 0.01 & 90 & 1822 & 100 & 0.01 & 63 & 3316 & 5 & 0.01 & 8 & 3244 & 5 & 0.01 & 8 \\ [-1.1ex]&0.99 & 877 & 100 & 0.01 & 94 & 1916 & 100 & 0.01 & 63 & 3780 & 4 & 0.01 & 7 & 3686 & 4 & 0.01 & 7 \\ \hline\multirow{6}{*}{D}&0.6 & 788 & 100 & 0.02 & 96 & 2114 & 100 & 0.01 & 60 & 4225 & 4 & 0.01 & 6 & 4145 & 4 & 0.01 & 5 \\ [-1.1ex]&0.7 & 673 & 100 & 0.02 & 97 & 2158 & 100 & 0.01 & 60 & 4302 & 4 & 0.01 & 5 & 4130 & 4 & 0.01 & 5 \\ [-1.1ex]&0.8 & 498 & 100 & 0.04 & 99 & 2210 & 100 & 0.01 & 63 & 4842 & 3 & 0.01 & 6 & 4634 & 3 & 0.01 & 6 \\ [-1.1ex]&0.9 & 328 & 100 & 0.2 & 99 & 2461 & 100 & 0.01 & 69 & 6319 & 3 & 0.01 & 7 & 5800 & 3 & 0.01 & 7 \\ [-1.1ex]&0.95 & 254 & 100 & 0.53 & 99 & 2208 & 100 & 0.01 & 82 & 9854 & 3 & 0.01 & 16 & 8816 & 3 & 0.01 & 12 \\ [-1.1ex]&0.99 & 186 & 100 & 3.13 & 99 & 5373 & 100 & 0.01 & 84 & 12413 & 5 & 0.01 & 32 & 9778 & 4 & 0.01 & 17 \\ \hline\multicolumn{2}{|c|}{All}&741 & 100 & 0.04 & 89 & 2105 & 100 & 0.01 & 63 & 4257 & 4 & 0.01 & 8 & 4088 & 4 & 0.01 & 7 \\\hline
\end{tabular}}
\caption{Pricing problem statistics of \texttt{CloudLarge} with 108 instances ($|\cN|=1000$)} \label{tab.price.CloudLarge}
\end{table}

 \bibliographystyle{elsarticle-num}
\bibliography{reference}


\end{document}